\documentclass[11pt,leqno]{article}

\usepackage{amsmath,amsfonts,amscd,amssymb,theorem}

\long\def\comment#1\endcomment{}

\makeatletter
\begingroup
\gdef\th@dotted{\normalfont\itshape
  \def\@begintheorem##1##2{%
        \item[\hskip\labelsep \theorem@headerfont ##1\ ##2.]}%
\def\@opargbegintheorem##1##2##3{%
   \item[\hskip\labelsep \theorem@headerfont ##1\ ##2\ (##3).]}}
\endgroup
\makeatother

\theoremstyle{dotted}

\newtheorem{theorem}{Theorem}[section]
\newtheorem{lemma}[theorem]{Lemma}

\newtheorem{prop}[theorem]{Proposition}
\newtheorem{corr}[theorem]{Corollary}

\makeatletter
\begingroup
\gdef\th@upshape{\normalfont
  \def\@begintheorem##1##2{%
        \item[\hskip\labelsep \theorem@headerfont ##1\ ##2.]}%
\def\@opargbegintheorem##1##2##3{%
   \item[\hskip\labelsep \theorem@headerfont ##1\ ##2\ (##3).]}}
\endgroup
\makeatother

\theoremstyle{upshape}

\newtheorem{defn}[theorem]{Definition}
\newtheorem{remark}[theorem]{Remark}
\newtheorem{exa}[theorem]{Example}

\makeatletter
\renewcommand{\subsection}{\@startsection{subsection}{2}{0pt}{-3ex
plus -1ex minus -0.2ex}{-2mm plus -0pt minus
-2pt}{\normalfont\bfseries}} 
\renewcommand{\subsubsection}{\@startsection{subsubsection}{3}{0pt}{-3ex
plus -1ex minus -0.2ex}{-2mm plus -0pt minus
-2pt}{\normalfont\bfseries}} 
\makeatother

\makeatletter
\@addtoreset{equation}{section}
\makeatother

\newcommand{\cntrct}                % contraction with a vector field
{\hspace{2pt}\raisebox{1pt}{\text{$\lrcorner$}}\hspace{2pt}}

\newcommand{\proof}[1][Proof.]{\smallskip\noindent{\em #1}}
\def\endproof{\hfill\ensuremath{\square}\par\medskip}

\renewcommand{\labelenumi}{{\normalfont(\roman{enumi})}}

\def\eqref#1{\thetag{\ref{#1}}}

\let\latexref=\ref
\def\ref#1{{\normalfont{\latexref{#1}}}}

\newcommand{\wt}{\widetilde}

\newcommand{\dg}{\dagger}

\newcommand{\idot}{{\:\raisebox{1pt}{\text{\circle*{1.5}}}}}
\newcommand{\hdot}{{\:\raisebox{3pt}{\text{\circle*{1.5}}}}}
\newcommand{\eps}{\varepsilon}
\renewcommand{\phi}{\varphi}

\newcommand{\bC}{\overline{C}}

\newcommand{\Hom}{\operatorname{Hom}}

\newcommand{\id}{\operatorname{\sf id}}
\newcommand{\Id}{\operatorname{\sf Id}}
\newcommand{\Iso}{\operatorname{\sf Iso}}

\newcommand{\gr}{\operatorname{\sf gr}}
\newcommand{\colim}{\operatorname{colim}}
\renewcommand{\lim}{\operatorname{lim}}

\newcommand{\D}{\mathcal{D}}
\newcommand{\F}{\mathcal{F}}
\newcommand{\C}{\mathcal{C}}

\newcommand{\K}{\mathcal{K}}

\newcommand{\Sets}{\operatorname{Sets}}
\newcommand{\End}{{\operatorname{End}}}

\newcommand{\amod}{{\text{\rm -mod}}}

\newcommand{\ppt}{{\sf pt}}

\newcommand{\lotimes}{\overset{\sf\scriptscriptstyle L}{\otimes}}

\newcommand{\Spec}{\operatorname{Spec}}
\newcommand{\QCoh}{\operatorname{QCoh}}
\newcommand{\Comm}{\operatorname{Comm}}
\newcommand{\Sch}{\operatorname{Sch}}

\newcommand{\Fun}{\operatorname{Fun}}

\newcommand{\E}{\mathcal{E}}

\newcommand{\copr}{\sqcup}

\newcommand{\Sec}{\operatorname{Sec}}

\newcommand{\Tw}{\operatorname{\sf tw}}
\newcommand{\Ar}{\operatorname{\sf ar}}

\newcommand{\Bi}{B_\infty}

\newcommand{\Cyl}{\operatorname{\sf C}}

\newcommand{\wCat}{\operatorname{Aug}}
\newcommand{\DwCat}{\operatorname{\mathcal{D}Aug}}
\newcommand{\Cat}{\operatorname{Cat}}
\newcommand{\iMor}{\Mor^*}
\newcommand{\Mor}{\operatorname{\mathcal{M}{\it or}}}
\newcommand{\DMor}{\operatorname{\mathcal{D}\mathcal{M}{\it or}}}
\newcommand{\cCat}{\operatorname{\mathcal{C}{\it at}}}

\newcommand{\FM}{\operatorname{\mathcal{F}\mathcal{M}}}
\newcommand{\FMbig}{\FM^{big}}

\newcommand{\bCat}{\overline{\Cat}}

\newcommand{\A}{\mathcal{A}}
\newcommand{\M}{\mathcal{M}}
\newcommand{\calo}{\mathcal{O}}

\newcommand{\rAdj}{\operatorname{Adj}}
\newcommand{\aAdj}{\operatorname{\mathcal{A}{\it dj}}}
\newcommand{\eEq}{\operatorname{\mathcal{E}{\it q}}}

\newcommand{\adj}{\operatorname{\sf adj}}
\newcommand{\Adj}{\operatorname{\sf Adj}}
\newcommand{\eq}{\operatorname{\sf eq}}
\newcommand{\Eq}{\operatorname{\sf Eq}}
\newcommand{\nat}{\operatorname{\sf nat}}
\newcommand{\Nat}{\operatorname{\sf Nat}}

\newcommand{\ev}{\operatorname{\sf ev}}

\newcommand{\ogamma}{\overline{\gamma}}
\newcommand{\wgamma}{\wt{\gamma}}

\newcommand{\Y}{\operatorname{\sf Y}}

\newcommand{\dm}{\diamond}

\newcommand{\R}{\operatorname{\sf R}}

\newcommand{\f}{\operatorname{\sf f}}
\newcommand{\aaa}{\operatorname{\sf a}}

\title{Adjunction in $2$-categories}

\author{D. Kaledin\thanks{Work supported by the Russian Science
    Foundation, grant 18-11-00141. No funding from any other
    source.}}

\begin{document}

\maketitle

\tableofcontents

\section*{Introduction.}

The fact that what categories form is a $2$-category rather than
simply a category is well-known but usually ignored. It is common to
see a commutative diagram of categories and functors, with the tacit
understanding that the diagram of course only commutes up to a fixed
isomorphism whose exact nature is obvious and left to the
reader. This is especially prominent when categories themselves
become an object of study, such as for instance in the recent
formulations of non-commutative algebraic geometry in the spirit of
\cite[Sections 2,3]{O}. Up to a point, this does not create any
problems. However, when the theory becomes more advanced, it can
become necessary to spell things out precisely. This is especially
the case when adjoint functors and various adjunction maps come into
play, such as for instance in the theory of spherical functors of
\cite{AL1}.

Unfortunately, spelling things out is not easy because the general
theory of $2$-categories stands on very cumbersome
foundations. Since everything only commutes up to an isomorphism,
these isomorphisms also have to be carried around, and one has to
constantly check the higher-order relations these isomorphisms are
supposed to satisfy. Even when one deals with monoidal categories
--- that is, $2$-categories with one object -- the resulting
combinatorics quickly gets completely out of hand (for a good
illustration, see e.g. Chapter VII of the classic textbook
\cite{Mcbook}). Giving a complete treatment of say adjunction in
such a language does not look like a thing an actual human can do.

Part of the reason for this situation could be the idea that
$2$-categories are inherently something that is one level up from
usual categories on a certain hypothetical ladder, whose higher
levels no-one ever saw, but that is supposed nevertheless to go all
the way to the equally hypothetic infinity. Thus theoretically, the
theory of $2$-categories is supposed to be built from scratch, and
its bizarre notation and arcane multilevel commutative diagrams are
not a bug but a feature: new combinatorics for a new world. But from
the practical viewpoint, this idea looks like a habitual illusion.
If one abandons it, then perfectly adequate tools for working with
$2$-categories are already available within the usual category
theory as known and practiced for decades. The first of these is
Grothendieck's packaging of ``families of categories'' into a
``Grothendieck fibration'' of \cite{sga}, so conspicuously and
inexplicably missing from \cite{Mcbook}. This is sometimes known as
the ``Grothendieck construction'', and at the end of the day, it is
one of the only two general constructions in abstract category
theory that are both non-trivial and useful for a working
mathematician (the other one being the notion of a factorization
system of \cite{bou}). The second ingredient is G. Segal's approach
to infinite loop spaces introduced in \cite{seg} that works for
categories just as well as, or even better than, for topological
spaces.

Combining the two ingredients, one obtains a simple and effective
way of packaging $2$-categories and related structures (monoidal
categories, symmetric monoidal categories) that has been around and
ready for use since at least early 1970-ies. In retrospect, it is a
pity that such a``Segal-style'' approach has not been adopted in the
foundational papers on Tannaka formalism such as \cite{De}. In any
case, for quite some time now, the approach has been completely
standard and ``well-known to the experts'' but not universally
known. Recently, with the advent of various homotopical enhancements
of category theory such as ``$\infty$-categories'', the situation
became both better and worse. On one hand, practitioners of that art
certainly only use the Segal approach since in that context, nothing
else could possibly work. On the other hand, among the uninitiated,
the perception now is that Segal-type approach to $2$-categories and
monoidal structures is something inextricably linked to fibrant
simplicial sets, quasicategories and the like, and the whole thing
comes as a take-it-or-leave-it monolith where everything depends on
everything else. The perception is utterly wrong but very common.

My own most recent encounter with the $2$-category formalism was in
the course of work on \cite{ka.hh} that expands and develops the
technology of trace theories and trace functors of \cite{trace},
with the goal of reaching very concrete applications such as
\cite[Section 3]{umn}. In \cite{ka.hh}, the Segal approach is simply
adopted by definition, without any explanations or
justifications. However, I feel that it might be useful to have a
write-up of how it is related to the more traditional
$2$-categorical technology. This is what the present paper tries to
provide.

The paper is mostly written in an overview style and contains
well-known results; the only thing that is possibly at least
partially new is the treatment of adjunction in
Section~\ref{adj.sec} that gave the name to the
paper. Section~\ref{cat.sec} mostly serves to fix notation and
terminology, and to give a brief introduction to the Grothendieck
construction with its main examples. Nothing whatsoever is new but
we have tried to make the exposition reasonably detailed and
self-contained. It would have been better to replace
Section~\ref{cat.sec} with a reference to a suitable textbook on
category theory but at present, none exist. We have also tried to be
very careful in spelling out various functorialities and
compatibility conditions, but we have allowed ourselves the liberty
of being quite cavalier about say homological algebra and derived
categories used in Example~\ref{sch.exa}, and various other things
very well served by existing literature. In similar style,
Section~\ref{2cat.sec} presents the Segal-style definitions of
various $2$-categorical notions and some related results. In
particular, we give a complete proof of the rigidification theorem
saying that any $2$-category is equivalent to a strict one (for the
precise statement, see Theorem~\ref{rigi.thm}). In
Section~\ref{adj.sec}, we treat adjunction. We do prove one result
skipped as obvious in \cite[Section 7]{ka.hh} (namely,
Lemma~\ref{adj.le}), but then we switch gears: for more technical
results in Section~\ref{adj.sec} that might also be new, we replace
proofs with precise references to \cite{ka.hh}. What we try to do is
to present things in a way that is maybe less canonical than
\cite[Section 7]{ka.hh} but more down-to-earth. Finally,
Section~\ref{app.sec} can be treated as one extended example: we
illustrate the general theory on one concrete real-life application,
that of the derived Morita $2$-category of DG algebras over a
commutative ring $k$, and the related Fourier-Mukai $2$-category of
algebraic varieties over $k$ as considered e.g.\ in \cite{AL2}. The
technology used is mostly that of \cite[Section 8]{ka.hh} (and it is
again not new but ``well-known to the experts''). Moreover, we again
allow ourselves to gloss over most of the homological details, and
we concentrate on checking that various canonical isomorphisms are
indeed canonical --- or rather, on explaining why in our formalism
there is almost nothing that needs to be actually checked. Since our
goal is to illustrate the general principle, we do not hesitate to
make simplifying assumptions --- thus Theorem~\ref{rigi.thm} is only
stated for small $2$-categories, and we only work with schemes of
finite type. The extended example, and the paper, ends with a
comparison theorem claiming that the Fourier-Mukai $2$-category over
some Noetherian commutative ring $k$ can be naturally embedded into
the derived Morita $2$-category over the same ring.

\subsection*{Acknowledgements.} I am very grateful to anonymous
referees and commentators on the Internet for many valuable
corrections and suggestions.

\section{Generalities}\label{cat.sec}

\subsection{Categories, functors and equivalences.}

Recall that by definition, a category $\C$ has a collection of
objects $c \in \C$, and a set $\C(c,c')$ of morphisms for any two
objects $c,c' \in \C$. Objects can form a set or a proper class; in
the former case, the category is called {\em small}. We also have a
distinguished identity morphism $\id = \id_c \in \C(c,c)$ for any $c
\in \C$, and composition maps $- \circ -:\C(c',c'') \times \C(c,c')
\to \C(c,c'')$ for any $c,c',c'' \in \C$ that are associative and
unital: we have
\begin{equation}\label{cat.eq}
(f \circ g) \circ h = f \circ (g \circ h), \qquad \id \circ f = f,
  \qquad f \circ \id = f
\end{equation}
for any composable triple of morphisms $f$, $g$, $h$. The {\em
  opposite category} $\C^o$ has the same objects, identities and
compositions as $\C$, but the morphisms are written in the opposite
direction: we have $\C^o(c,c') = \C(c',c)$. For any morphism $f \in
\C(c,c')$, we let $f^o$ be the same $f$ but considered as a morphism
in the opposite category. An object $o \in \C$ is {\em initial} if
$\C(o,c) \cong \ppt$ for any $c \in \C$. If an initial object
exists, then it is unique up to a unique isomorphism, and we can add
a new initial object $o$ to any category $\C$. We denote the
resulting category by $\C^<$. Dually, a terminal object $o \in \C$
is an initial object of $\C^o$, and we can always add a new terminal
object to any $\C$ and obtain a category $\C^> = (\C^{o<})^o$.

A functor $\gamma:\C \to \C'$ associates an object $\gamma(c) \in
\C'$ to any object $c \in \C$, and a map $\C(c,c') \to
\C'(\gamma(c),\gamma(c'))$ to any two objects $c,c' \in \C$; the
maps commute with compositions and send identity morphisms to
identity morphisms. A functor $\gamma:\C_0 \to \C_1$ induces the
opposite functor $\gamma^o:\C_0^o \to \C_1^o$ between the opposite
categories. For any category $\C$, we have natural functors $\C \to
\C^<$, $\C \to \C^>$.

Sets and maps between them form a category denoted $\Sets$. Small
categories and functors between them form a category denoted
$\Cat$. A category is discrete if all its morphisms are identity
maps. A discrete small category is the same thing as a set, and
sending a set $S \in \Sets$ to itself considered as a discrete
category gives an embedding $\Sets \subset \Cat$. Alternatively, one
can consider the category whose objects are elements $s \in S$, and
with exactly one morphism between any two objects. We will denote
this category by $e(S)$. Any partially ordered set $J$ can be
considered as a category with $J(j,j')=\ppt$ if $j' \geq j$ and
$J(j,j') = \emptyset$ otherwise. In particular, for any integer $n
\geq 0$, we have the totally ordered set $[n] = \{0,\dots,n\}$, with
the usual order. If $n=1$, then $[1]$ is the ``single-arrow
category'' with two objects $0$, $1$ and unique non-identity
morphism $0 \to 1$.

Recall that a subcategory $\C' \subset \C$ is {\em full} if
$\C'(c,c') = \C(c,c')$ for any $c,c' \in \C'$ -- that is, $\C$ and
$\C'$ have the same morphisms -- and we will say that $\C' \subset
\C$ is {\em dense} if they have the same objects. Giving a dense
subcategory is equivalent to giving a class $v$ of maps in $\C$ that
is closed under compositions and contains all the identity maps. For
brevity, we call such classes {\em closed}, and we denote the
corresponding dense subcategory by $\C_v \subset \C$. The minimal
closed class is the class $\Id$ consisting only of identity maps,
and $\C_{\Id} \subset \C$ is the discrete category underlying
$\C$. Another useful closed class is the class $\Iso$ of all
invertible maps.

A morphism, or a map, or a natural transformation $\alpha$ between
two functors $\gamma_0,\gamma_1:\C \to \C'$ is a collection of maps
$\alpha(c):\gamma_0(c) \to \gamma_1(c)$, $c \in \C$, such that
$\alpha(c') \circ \gamma_0(f) = \gamma_1(f) \circ \alpha(c)$ for any
map $f:c \to c'$ in $\C$. The same data can be packaged into a
single functor
\begin{equation}\label{ga.1}
\gamma:[1] \times \C \to \C'
\end{equation}
whose restriction to $0 \times \C$ resp.\ $1 \times \C$ is
identified with $\gamma_0$ resp.\ $\gamma_1$. A map $\alpha$ is an
isomorphism if so is $\alpha(c)$ for any $c \in \C$, and this
happens if and only if the functor \eqref{ga.1} extends to
$e(\{0,1\}) \times \C \supset [1] \times \C$.

For any two categories $\C_0$, $\C_1$ and functors $\pi_0:\C_0 \to
\C$, $\pi_1:\C_1 \to \C$ to a third category $\C$, we will denote by
$\C_0 \times_\C \C_1$ the category of triples $\langle
c_0,c_1,\alpha \rangle$ of objects $c_0 \in \C_0$, $c_1 \in \C_1$
and an isomorphism $\alpha:\pi_0(c_0) \cong \pi_1(c_1)$. Note that
this notation is somewhat abusive, since even when $\C_0$, $\C_1$,
$\C$ are small, $\C_0 \times_\C \C_1$ is {\em not} the fibered
product in $\Cat$: the latter is the subcategory of triples $\langle
c_0,c_1,a \rangle$ such that $\alpha \in \Id$. Analogously, a {\em
  functor $\C_0 \to \C_1$ over $\C$} is a pair $\langle
\gamma,\alpha \rangle$ of a functor $\gamma:\C_0 \to \C_1$ and an
isomorphism $\alpha:\pi_1 \circ \gamma \cong \pi_0$, and for any $\C
\in \Cat$, we denote by $\Cat/\C$ the category whose objects are
small categories $\C'$ equipped with a functor $\C' \to \C$, and
whose morphisms are functors over $\C$. A functor $\langle
\gamma,\alpha \rangle:\C_0 \to \C_1$ over $\C$ is {\em strict} if
$\alpha(c) \in \Id$ for any $c \in \C_0$ (that is, $\gamma$ strictly
commutes with projections to $\C$). A functor $\gamma:\C_0 \to \C$
induces the pullback functor $\gamma^*:\Cat/\C \to \Cat/\C_0$
sending $\C' \to \C$ to $\gamma^*\C' = \C_0 \times_\C \C'$. The {\em
  fiber} $\C'_c$ of a functor $\C' \to \C$ over an object $c \in \C$
is the subcategory $\C'_c \subset \C$ spanned by objects $c'$ such
that $\pi(c') = c$ and morphisms $f$ such that $\pi(f) = \id_c$.

A functor $\gamma:\C \to \C'$ is an equivalence if there exists a
functor $\gamma':\C' \to \C$ and isomorphisms $a:\Id \cong \gamma'
\circ \gamma$, $a':\gamma \circ \gamma' \cong \Id$. This is a
condition and not an extra structure, and the triple $\langle
\gamma',a,a' \rangle$ can be chosen in a canonical way: if some such
triple exists, then $\gamma$ admits a right-adjoint functor
$\gamma^\dg$, and the adjunction maps $\Id \to \gamma^\dg \circ
\gamma$, $\gamma \circ \gamma^\dg \to \Id$ are isomorphisms. We
recall that the adjoint functor, if it exists, is unique up to a
unique isomorphism. Alternatively, one can take the left-adjoint
functor, and the isomorphisms inverse to the adjunction maps.

A functor is an equivalence if and only if it is fully faithful and
essentially surjective (that is, surjective on isomorphism classes
of objects). The class of all equivalences is closed in $\Cat$, and
in practice, one is only interested in categories ``up to an
equivalence''. To make this formal, one can invert all equivalences
in $\Cat$. Namely, say that a functor $\gamma:\C \to \C'$ inverts a
morphism $f$ in $\C$ if $\gamma(f)$ is invertible in $\C$', and
denote by $\bCat$ the category of small categories and isomorphism
classes of functors between them. We then have the tautological
functor $\tau:\Cat \to \bCat$ sending a category to itself and a
functor to its isomorphism class, and by definition, this functor
inverts all equivalences in $\Cat$.

\begin{lemma}\label{bcat.le}
Any functor $\gamma:\Cat \to \C$ to some category $\C$ that inverts
all equivalences factors as
$$
\begin{CD}
\Cat @>{\tau}>> \bCat @>{\gamma'}>> \C,
\end{CD}
$$
and the factorization is unique up to a unique isomorphism.
\end{lemma}

\proof{} To prove existence, we can take $\gamma'=\gamma$ on
objects; we then need to check that $\gamma(F_0) = \gamma(F_1)$ as
soon as two functors $F_0,F_1:I \to I'$ between two categories $I,I'
\in \Cat$ are isomorphic. Indeed, such a pair gives a functor
$F:e(\{0,1\}) \times I \to I'$, and the embeddings $i_l:I \to
e(\{0,1\}) \times I$ onto $l \times I$, $l=0,1$ are equivalences
with common right-inverse given by the projection $e(\{0,1\}) \times
I \to I$. Therefore $\gamma(i_0) = \gamma(i_1)$, so that
$\gamma(F_0) = \gamma(F_1)$. Uniqueness is left to the reader.
\endproof

\subsection{Strict $2$-categories and $2$-functors.}

A strict $2$-category $\C$ is a ``category enriched in $\Cat$'': it
has objects $c \in \C$, a category of morphisms $\C(c,c')$ for any
two objects $c,c' \in \C$, identity objects $\id = \id_c \in
\C(c,c)$, and composition functors $\C(c',c'') \times \C(c,c') \to
\C(c,c'')$ such that \eqref{cat.eq} holds on the nose. To avoid
set-theoretical difficulties, it is prudent to assume that all the
categories $\C(c,c')$ are small. If objects form a set and not a
proper class, then $\C$ itself is called small.

For a strict $2$-category $\C$, there are two opposite
$2$-categories that we call the {\em $1$-opposite} $2$-category
$\C^\iota$ and the {\em $2$-opposite} $2$-category $\C^\tau$. They
all have the same objects, and morphism categories
$\C^\iota(c,c')=\C(c',c)$, $\C^\tau(c,c') = \C(c,c')^o$.

A usual category $\C$ is trivially a strict $2$-category, with discrete
$\C(c,c')$. The first non-trivial example is the $2$-category
$\cCat$ of small categories: for any two small categories $I$, $I'$,
$\cCat(I,I')$ is the category $\Fun(I,I')$ of functors from $I$ to
$I'$. More generally, for any $I \in \cCat$, functors over $I$ also
form a category in the obvious way, and we have the $2$-category
$\cCat/I$: objects are small categories $\C$ equipped with a functor
$\C \to I$, and $\cCat/I(\C,\C')$ is the category $\Fun_I(\C,\C')$
of functors from $\C$ to $\C'$ over $I$.

Strict $2$-functors are defined in a straightforward manner, but for
most practical applications, they are too strict. It is more
convenient to use the following notion.

\begin{defn}\label{st.2fun.def}
Assume given strict $2$-categories $\C$ and $\C'$. A {\em $2$-functor}
$\gamma:\C \to \C'$ is a rule that associates an object $\gamma(c)
\in \C'$ to any object $c \in \C$, a functor $\gamma_{c,c'}:\C(c,c')
\to \C'(\gamma(c),\gamma(c'))$ to any pair of objects $c,c' \in \C$,
an isomorphism $\gamma_c:\id_{\gamma(c)} \cong \gamma_{c,c}(\id_c)$
to any object $c \in \C$, and an isomorphism $\gamma_{c,c',c''}:\mu'
\circ (\gamma_{c',c''} \times \gamma_{c,c'}) \cong \gamma_{c,c''}
\circ \mu$ to any triple of objects $c,c',c'' \in \C$, where $\mu$
resp.\ $\mu'$ stands for the composition functors in $\C$
resp.\ $\C'$. For any objects $c,c' \in \C$ and morphism $f \in
\C(c,c')$, we must have
$$
\gamma_{c,c,c'}(f \times \id_c) \circ (\id \circ \gamma_c) = \id,
\qquad \gamma_{c,c',c'}(\id_{c'} \times f) \circ (\gamma_{c'} \circ
\id) = \id,
$$
and for any four objects $c,c',c'',c''' \in \C$ and three morphisms
$f \in \C(c,c')$, $f' \in \C(c',c'')$, $f'' \in \C(c'',c''')$, the
square
$$
\begin{CD}
\gamma_{c'',c'''}(f'') \circ \gamma_{c',c''}(f') \circ
\gamma_{c,c'}(f) @>{\gamma_{c',c'',c'''} \circ \id}>>
\gamma(c',c''')(f'' \circ f') \circ\gamma_{c,c'}(f)\\
@V{\id \circ \gamma_{c,c',c''}}VV @VV{\gamma_{c,c',c'''}}V\\
\gamma_{c'',c'''}(f'') \circ \gamma_{c,c''}(f' \circ f)
@>{\gamma_{c,c'',c'''}}>> \gamma_{c,c'''}(f'' \circ f' \circ f)
\end{CD}
$$
must be commutative.
\end{defn}

The collection of $2$-functors $\gamma:\C \to \C'$ between strict
$2$-categories $\C$, $\C'$ carries a rich and varied structure: one
can define morphisms $\gamma \to \gamma'$ by using \eqref{ga.1},
then morphisms between morphisms, and then all sorts of compositions
and associativity constraints. We will need only a small part of
this cornucopia, namely, the following.

\begin{defn}\label{fun2.def}
Assume given $2$-functors $\gamma,\gamma':\C \to \C'$ between strict
$2$-categories such that $\gamma(c)=\gamma'(c)$ for any $c \in
\C$. A {\em $2$-morphism} $\alpha:\gamma \to \gamma'$ is a
collection of morphisms $\alpha_{c,c'}:\gamma_{c,c'} \to
\gamma'_{c,c'}$, $c,c' \in \C$ such that for any triple of objects
$c,c',c'' \in \C$, the square
$$
\begin{CD}
\mu' \circ (\gamma_{c,c'} \times \gamma_{c',c''})
@>{\gamma_{c,c',c''}}>> \gamma_{c,c''} \circ \mu\\
@V{\mu' \circ (\alpha_{c,c'} \times \alpha_{c',c''})}VV @VV{\alpha_{c,c''}}V\\
\mu' \circ (\gamma'_{c,c'} \times \gamma'_{c',c''})
@>{\gamma'_{c,c',c''}}>> \gamma'_{c,c''} \circ \mu
\end{CD}
$$
is commutative. A $2$-morphism $\alpha$ is {\em invertible} if so
are all the maps $\alpha_{c,c'}$.
\end{defn}

With this definition, as soon as $\C$ is small, $2$-functors $\C \to
\C'$ and $2$-morphisms between them form a category that we denote
$\Fun^2(\C,\C')$ (if $2$-functors $\gamma,\gamma':\C \to \C'$ do not
coincide on objects, there are no maps $\gamma \to \gamma'$ in
$\Fun^2(\C,\C')$). For any strict $2$-category $\C''$,
postcomposition with a $2$-functor $\C' \to \C''$ defines a functor
$\Fun^2(\C,\C') \to \Fun^2(\C,\C'')$, and if $\C''$ is small, then
precomposition with a $2$-functor $\gamma:\C'' \to \C$ defines a
functor $\gamma^*:\Fun^2(\C,\C') \to \Fun^2(\C'',\C')$. Note that if
$\C$ and $\C'$ are usual categories considered as $2$-categories
with discrete categories of morphisms, then $\Fun^2(\C,\C')$ is
discrete: objects are functors, but we do not allow any non-trivial
maps between them.

\begin{defn}\label{stro.def}
A $2$-functor $\gamma:\C \to \C'$ is a {\em strong equivalence} if
there exists a $2$-functor $\gamma':\C' \to \C$ and invertible
$2$-morphisms $\id \to \gamma' \circ \gamma$, $\gamma \circ \gamma'
\to \id$.
\end{defn}

We note that strong equivalences in the sense of
Definition~\ref{stro.def} are indeed rather strong -- in particular,
objects in $\C$ and $\C'$ are canonically the same (and identified
by $\gamma$). For usual small categories considered as
$2$-categories, strong equivalence gives an isomorphism in $\Cat$,
not just an equivalence.

\subsection{Grothendieck construction.}\label{groth.subs}

Definition~\ref{st.2fun.def} is a generalization of the notion of a
{\em pseudofunctor} introduced by Grothendieck in \cite{sga}. A
pseudofunctor is a $2$-functor from a usual category $I$ to the
$2$-category $\cCat$ of small categories. The great discovery of
\cite{sga} is that all the data defining a pseudofunctor can be very
efficiently packaged by means of the usual category theory. These
day, it is known as the {\em Grothendieck construction}. Here is a
brief summary.

\begin{defn}\label{fib.def}
A morphism $f:c' \to c$ in a category $\C$ is {\em cartesian} with
respect to a functor $\pi:\C \to I$ if any $f':c'' \to c$ such that
$\pi(f) = \pi(f')$ uniquely factors as $f' = f \circ f_0$ with
$\pi(f_0)=\id$. A {\em cartesian lifting} of a morphism $f:i' \to i$
in $I$ is a morphism $f'$ in $\C$ cartesian with respect to $\pi$
and such that $\pi(f')=f$. A functor $\pi:\C \to I$ is a {\em
  prefibration} if for any $c \in \C$, any morphism $f:i' \to
i=\pi(c)$ in $I$ admits a cartesian lifting $f':c' \to c$. A
prefibration is a {\em fibration} if the composition of cartesian
morphisms is cartesian. A morphism is {\em cocartesian} if it is
cartesian as a morphism in the opposite category $\C^o$ with respect
to the opposite functor $\pi^o:\C^o \to I^o$, and a functor $\pi$ is
a {\em cofibration} if the opposite functor $\pi^o$ is a fibration.
\end{defn}

\begin{remark}\label{comma.rem}
For any functor $\pi:\C \to I$ and object $i \in I$, the {\em right
  comma-fiber} $i \setminus^\pi \C$ is the category of pairs
$\langle c,\alpha \rangle$ of an object $c \in \C$ and a map
$\alpha:i \to \pi(c)$. Then by definition, $\C_i \subset i
\setminus^\pi \C$ is the full subcategory spanned by pairs $\langle
c,\alpha \rangle$ with $\alpha = \id_i$, and $\pi$ is a prefibration
if and only for any $i \in I$, the embedding $\C_i \to i
\setminus^\pi \C$ admits a right-adjoint. In particular, if $i \in
I$ is an initial object, then $i \setminus^\pi \C \cong \C$, and the
embedding $\C_i \subset \C$ admits a right-adjoint. Dually, one has
a charaterization of cofibrations in terms of left comma-fibers
$\C/^\pi i = (i \setminus^\pi \C^o)^o$, and if $i \in I$ is
terminal, then for any cofibration $\C \to I$, the embedding $\C_i
\subset \C$ admits a left-adjoint.
\end{remark}

While cartesian liftings in Definition~\ref{fib.def} are not unique,
they are unique up to a unique isomorphism, by virtue of the
universal property of cartesian maps. Therefore for any prefibration
$\pi:\C \to I$ and morphism $f:i' \to i$ in $I$, sending $c \in \C_i$ to
the source $c'$ of the cartesian lifting $f:c' \to c$ of the map $f$
defines a functor $f^*:\C_{i'} \to \C_i$ known as the {\em transition
  functor} of the prefibration. By the same universal property, we
have functorial maps
\begin{equation}\label{fib.eq}
\Id \to \id_i^*, \qquad g^* \circ f^* \to (f \circ g)^*
\end{equation}
for any $i \in I$ and any composable pair of maps $f$, $g$ in
$I$. The first of these maps is always an isomorphism, and the
second is an isomorphism for any $f$, $g$ if and only if $\pi$ is a
fibration. In the latter case, the correspondences $i \mapsto \C_i$,
$f \mapsto f^*$ and the maps \eqref{fib.eq} define a $2$-functor
$\gamma:I^o \to \cCat$, with all the conditions of
Definition~\ref{st.2fun.def} again satisfied by virtue of the
universal property of cartesian maps. Conversely, for any
$2$-functor $\gamma:I^o \to \cCat$, we can define a category $\C$
whose objects are pairs $\langle i,c \rangle$, $i \in I$, $c \in
\C_i = \gamma(i)$, and whose morphisms from $\langle i,c \rangle$ to
$\langle i',c' \rangle$ are pairs of a map $f:i \to i'$ and a map $c
\to \gamma(f^o)(c')$. Then the forgetful functor $\C \to I$,
$\langle i,c \rangle \mapsto i$ is a fibration with transition
functors $\gamma(f^o)$. The two constructions are mutually inverse,
so they describe the same entity, but the fibration description is
much more economical: the only data are the category $\C$ and the
functor $\C \to I$. One does not have to keep track of all the
isomorphisms and compatibility conditions of
Definition~\ref{st.2fun.def}.

For a cofibration $\C \to I$, the picture is dual: we have a
transition functor $f_!:\C_i \to \C_{i'}$ for any map $f:i \to i'$
in $I$, and we have the corresponding versions of the isomorphisms
\eqref{fib.eq} (both go in the opposite direction, but since they
are invertible, this is not so important).

\begin{exa}\label{cyl.exa}
Assume given two categories $\C_0$, $\C_1$, and a functor
$\gamma:\C_0 \to \C_1$. Then the {\em cylinder} $\Cyl(\gamma)$ is
the category whose objects are those of $\C_0$ and of $\C_1$, and
whose morphisms are given by
$$
\Cyl(\gamma)(c,c')=\begin{cases}
\C_l(c,c'), &\quad c,c' \in \C_l, l=0,1\\
\C_1(\gamma(c),c'), &\quad c \in \C_0,c' \in \C_1,
\end{cases}
$$
and no morphisms from $c \in \C_1$ to $c' \in \C_0$. We have the
natural projection $\Cyl(\gamma) \to [1]$ with fibers $\C_0$,
$\C_1$, it is a cofibration with transition functor $\gamma$, and
every cofibration $\C \to [1]$ is of this form.
\end{exa}

By definition, a cofibration $\pi:\C \to I$ defines the opposite
fibration $\pi^o:\C^o \to I^o$. It also defines the {\em transpose
  fibration} $\C^\perp \to I^o$ with fibers $\C_i$, transition
functors $f^* = f^o_!$, and the isomorphisms \eqref{fib.eq} inverse
to the corresponding isomorphisms for $\pi$. Moreover, it might
happen that $\pi$ itself is a fibration (it is then called a {\em
  bifibration}). This happens if and only if all the cofibration
transition functors $f_!$ admit right-adjoint functors $f^*$, and
these right-adjoints are then the fibration transition
functors. Dually, a fibration $\C \to I$ might also be a
cofibration, and in any case, it defines the opposite cofibration
$\C^o \to I^o$ and the transpose cofibration $\C_\perp \to I^o$. One
can also combine the constructions: for any fibration $\C \to I$,
the opposite fibration $\C_\perp^o \to I$ to the transpose
cofibration $\C_\perp \to I^o$ has fibers $(\C_\perp^o)_i = \C_i^o$,
and transition functors $(f^*)^o$, and similarly for cofibrations.

\begin{remark}
The definitions of the transpose fibrations and cofibrations given
above use the Grothendieck construction twice, and this is
unpleasant since it involves choices. Although choices do cancel out
in the end, here is a more direct definition. For any fibration
$\pi:\C \to I$, the transpose cofibration $\C_\perp$ has the same
objects, and morphisms are isomorphism classes of diagrams
\begin{equation}\label{dom.0}
\begin{CD}
c @<{f}<< \wt{c} @>{v}>> c',
\end{CD}
\end{equation}
where $f$ is cartesian with respect to $\pi$, and $\pi(v)$ is
invertible. Compositions are obtained by taking fibered products,
and the existence of the relevant fibered products easily follows
from the definition of a fibration. For the transpose fibration
$\C^\perp$ to a cofibration $\C \to I$, the arrows in \eqref{dom.0}
go in the other direction, and one uses coproducts rather than
products for the compositions.
\end{remark}

\subsection{Cartesian functors.}

If we have two fibrations $\C,\C' \to I$, then a functor $\gamma:\C
\to \C'$ over $I$ is {\em cartesian over a map} $f:i \to i'$ if it
sends all cartesian liftings of $f$ to maps in $\C'$ cartesian over
$I$, and {\em cartesian} if it is cartesian over all maps (that is,
sends cartesian maps in $\C$ to cartesian maps in $\C'$). To see
this condition more explicitly, it is useful to observe that more
generally, if we have two functors $\C,\C' \to I$, and $\C' \to I$
is a fibration resp.\ a cofibration, then a functor $\langle
\gamma,\alpha \rangle$ from $\C$ to $\C'$ over $I$ is isomorphic to
a strict functor over $I$, unique up to a unique isomorphism and
called the {\em strictification} of $\gamma$ (this happens because
just like any other map, the isomorphisms $\alpha(c)$, $c \in \C$
admit cartesian resp.\ cocartesian liftings to $\C'$). If $\C,\C'
\to I$ are fibrations, then up to a unique isomorphism, a functor
$\gamma:\C \to \C'$ over $I$ is expicitly given by a collection of
functors $\gamma(i):\C_i \to \C'_i$, $i \in I$, induced by its
strictification, and functorial maps
\begin{equation}\label{fib.fu}
\gamma(f):\gamma(i) \circ f^* \to f^* \circ \gamma(i')
\end{equation}
for all morphisms $f:i \to i'$ in $I$, subject to natural
compatibility conditions that we leave to the reader. The functor
$\gamma$ is cartesian over $f$ if \eqref{fib.fu} is an
isomorphism. If $\C$ is small, then cartesian functors form a full
subcategory $\Fun^\natural_I(\C,\C') \subset \Fun_I(\C,\C')$ in the
category of all functors over $I$, and if we restrict our attension
to strict functors and/or strict cartesian functors, we obtain
equivalent categories. Fibrations and strict cartesian functors then
form a subcategory $\Cat^\natural/I \subset \Cat/I$ and a
$2$-subcategory $\cCat^\natural/I \subset \cCat/I$, and the
$2$-subcategory spanned by all cartesian functors is strongly
equivalent to $\cCat^\natural/I$.

Dually, for cofibrations, a functor $\gamma:\C \to \C'$ induces
functors $\gamma(i)$, $i \in I$, and maps $f_! \circ \gamma(i') \to
\gamma(i) \circ f_!$ dual to \eqref{fib.fu}, unique up to a unique
isomorphism, a functor is cocartesian iff all these maps are
isomorphisms, and cofibrations together with strict cocartesian
functors form a subcategory $\Cat_\natural/I \subset \Cat/I$ and a
$2$-subcategory $\cCat_\natural/I \subset \cCat/I$, where in the
latter case, dropping ``strict'' gives a strongly equivalent
$2$-category. The transpose fibration/transpose cofibration
construction provides an equivalence
\begin{equation}\label{fi.co.eq}
\Cat^\natural/I \cong \Cat_\natural/I^o.
\end{equation}
For any fibration $\C \to I$ and any functor $\gamma:I' \to I$, the
pullback $\gamma^*\C \to I'$ is a fibration, a pullback of a functor
cartesian over $I$ is cartesian over $I'$, and similarly for
cofibrations and cocartesian functors.

\begin{remark}\label{iso.rem}
As another instance of the strictification phenomenon for fibrations
and cofibrations, we note that up to a canonical equivalence, the
pullback $\gamma^*\C$ of a fibration or a cofibration of small
categories is given by an honest fibered product in $\Cat$.
\end{remark}

\begin{exa}\label{adj.perp.exa}
As a motivation for the transpose cofibration construction, assume
given two fibrations $\C,\C' \to I$, and a functor $\gamma:\C \to
\C'$ over $I$ such that for any $i \in I$, the component
$\gamma(i):\C_i \to \C'_i$ admits a right-adjoint
$\gamma^\dg(i):\C'_i \to \C_i$. Then for any map $f:i \to i'$ in
$I$, the map \eqref{fib.fu} induces by adjunction a map $f^* \circ
\gamma^\dg(i') \to \gamma^\dg(i) \circ f^*$, but these maps do {\em
  not} define a functor from $\C'$ to $\C$ over $I$: they go in the
wrong direction. What they define is a functor
$\gamma^\dg_\perp:\C'_\perp \to \C_\perp$ over $I^o$ between the
transpose cofibrations. If $\gamma^\dg_\perp$ is cocartesian, then
it is transpose to a cartesian functor $\gamma^\dg:\C' \to \C$
right-adjoint to $\gamma$. If $\gamma$ is cartesian, then its
transpose functor $\gamma_\perp:\C_\perp \to \C'_\perp$ is
left-adjoint to $\gamma^\dg_\perp$. But in general, while
$\gamma^\dg_\perp$ is defined by $\gamma$ uniquely up to a unique
isomorphism, it is not adjoint to $\gamma$ in any obvious sense.
\end{exa}

A {\em section} of a functor $\C \to I$ is a functor $I \to \C$ over
$I$, and if $\C \to I$ is a fibration, a section is cartesian if it sends
all maps in $I$ to cartesian maps. If $I$ is small, then sections
form a category that we denote by $\Sec(I,\C)$, and if $\C \to I$ is
a fibration, we denote by $\Sec^\natural(I,\C) \subset \Sec(I,\C)$
the full subcategory spanned by cartesian sections. For any
subcategory $I' \subset I$, with the embeddng functor $j:I' \to I$,
we have the restriction functor
\begin{equation}\label{sec.I}
\Sec^\natural(I,\C) \to \Sec^\natural(I',j^*\C).
\end{equation}
A useful observation is that if $I' \subset I$ is full, and the
embedding $j$ admits a left-adjoint functor $j^\dg:I \to I'$, then
\eqref{sec.I} is an equivalence. To see this, we may restrict our
attention to strict sections. Then the inverse equivalence sends a
strict section $\sigma \in \Sec^\natural(I',j^*\C)$ to the strict
section $\sigma' \in \Sec(I,\C)$ given by $\sigma'(i) =
a^*\sigma(j^\dg(i))$, $i \in I$, where $a:i \to j(j^\dg(i))$ is the
adjunction map, with the obvious maps \eqref{fib.fu}. Dually, if $\C
\to I$ is a cofibration, we let $\Sec_\natural(I,\C) \subset
\Sec(I,\C)$ be the full subcategory spanned by cocartesian sections,
and we have the equivalence \eqref{sec.I} whenever the full
embedding $I' \to I$ admits a right-adjoint.

A composition $\pi' = \pi \circ \gamma:\C' \to I$ of two fibrations
$\gamma:\C' \to \C$, $\pi:\C \to I$ is a fibration but the converse
is not true: $\pi'$ can be a fibration even if $\gamma$ is not. One
has to impose additional conditions; here is a list.

\begin{lemma}\label{fib.le}
Assume given a fibration $\pi:\C \to I$ and a functor $\gamma:\C'
\to \C$. Then $\gamma$ is a fibration if and only if \thetag{i}
$\pi' = \pi \circ \gamma:\C' \to I$ is a fibration, \thetag{ii}
$\gamma$ is cartesian over $I$, \thetag{iii} for any $i \in I$,
$\gamma(i)$ is a fibration, and \thetag{iv} for any map $f:i \to i'$
in $I$, the functor $\C'_{i'} \to (f^*)^*\C'_i$ induced by the
commutative diagram
\begin{equation}\label{fib.dia}
\begin{CD}
\C'_{i'} @>{f^*}>> \C'_i\\
@V{\gamma(i')}VV @VV{\gamma(i)}V\\
\C_{i'} @>{f^*}>> \C_i
\end{CD}
\end{equation}
is cartesian over $\C_{i'}$.
\end{lemma}

\proof{} The ``only if'' part is easy and left to the reader. For
the ``if'' part, note that each condition makes sense only if the
previous ones are satisfied (in particular, \eqref{fib.dia} in
\thetag{iv} is commutative by virtue of \thetag{ii}). If we let
$C$ be the class of maps in $\C$ cartesian over $I$, and let $V$ be
the class of maps $f$ that are ``vertical'' in the sense that
$\pi(v)$ is invertible, then \thetag{i} and \thetag{ii} insure that
$\gamma$ is a fibration over $\C_C$, and \thetag{iii} insures that
it is a fibration over $\C_V$. By the definition of a fibration,
every map $f$ in $\C$ factors as $f = c \circ v$, $v \in V$, $c \in
C$, uniquely up to a unique isomorphism, and \thetag{i},
\thetag{ii}, \thetag{iii} together insure that $\gamma$ is a
prefibration with transition functors $f^* \cong v^* \circ
c^*$. Moreover, the second map in \eqref{fib.eq} is an isomorphism
if $f \in C$, $g \in C$, or $f \in C$, $g \in V$, or $f \in V$, $g
\in V$. It then remains to check that it is an isomorphism when $f
\in V$ and $g \in C$, and this immediately follows from \thetag{iv}.
\endproof

\subsection{Examples.}

Let us now give some general examples of fibrations and cofibrations
in the sense of Definition~\ref{fib.def}.

\begin{exa}\label{ar.exa}
For any category $I$, the {\em arrow category} $\Ar(I) =
\Fun([1],I)$ is the category of arrows $i \to i'$ in $I$, with
morphisms from $f_0:i_0 \to i_0'$ to $f_1:i_1 \to i_1'$
given by commutative squares
\begin{equation}\label{ar.sq}
\begin{CD}
i_0 @>{f_0}>> i_0'\\
@V{g}VV @VV{g'}V\\
i_1 @>{f_1}>> i_1'.
\end{CD}
\end{equation}
Consider the functors $s,t:\Ar(I) \to I$ sending an arrow to its
source resp.\ its target. Then $s$ is a fibration, its fibers
$\Ar(I)_i$ are the right-comma fibers $i \setminus^{\Id} I$ of
Remark~\ref{comma.rem} of the identity functor $\Id:I \to I$, and
and a morphism $f$ is cartesian with respect to $s$ if and only if
$t(f)$ is invertible. Dually, $t$ is a cofibration with fibers
$\Ar(I)_i \cong I /^{\Id} i$, and $f$ is cocartesian if and only if
$s(f)$ is invertible.
\end{exa}

\begin{exa}\label{sq.exa}
For any class of morphisms $V$ in $I$, let $\Ar^V(I) \subset \Ar(I)$
be the full subcategory spanned by arrows in $V$. Say that $V$ is
{\em closed under pullbacks} if for any $f_1:i_1 \to i_1'$ in $V$
and any map $g':i_0' \to i_1'$, there exists a cartesian square
\eqref{ar.sq} in $I$ with $f_0 \in V$. Then in this case, the
projection $t:\Ar^V(I) \to I$ of Example~\ref{ar.exa} is a
fibration, with cartesian maps corresponding to cartesian squares
\eqref{ar.sq}.
\end{exa}

\begin{exa}\label{tw.exa}
For any category $I$, the fibration $t^\perp:\Tw(I) = \Ar(I)^\perp
\to I^o$ transpose to the cofibration $t:\Ar(I) \to I$ is the {\em
  twisted arrow category} of the category $I$: its objects are
arrows $f:i \to i'$, and morphisms from $f_0:i_0 \to i_0'$ to $f_1:i_1
\to i_1'$ are given by commutative diagrams
$$
\begin{CD}
i_0 @>{f_0}>> i'_0\\
@VVV @AAA\\
i_1 @>{f_1}>> i'_1.
\end{CD}
$$
Note that the projection $s \times t^\perp:\Tw(I) \to I \times I^o$
is also a fibration; it corresponds to the functor $I(-,-):I^o
\times I = (I \times I^o)^o \to \Sets \subset \Cat$ by the
Grothendieck construction.
\end{exa}

\begin{exa}\label{facto.exa}
By definition (see e.g. \cite{bou}), a {\em factorization system}
$\langle L,R \rangle$ on a category $I$ consists of two closed
classes of maps $L$, $R$ in $I$ such that $L \cap R = \Iso$, and any
morphism $f:i \to i'$ factors as
\begin{equation}\label{facto.dia}
\begin{CD}
i @>{l}>> \wt{i} @>{r}>> i',
\end{CD}
\end{equation}
with $l \in L$, $r \in R$, uniquely up to a unique isomorphism. Then
for any factorization system $\langle L,R \rangle$, the projection
$s:\Ar^L(I) \to I$ of Example~\ref{ar.exa} is a fibration, and a map
$f$ is cartesian iff $t(f) \in R$. Dually, $t:\Ar^R(I) \to I$ is a
cofibration, and a map $f$ is cocartesian iff $s(f) \in L$.

Moreover, if we are given a fibration $\pi:I' \to I$, then the
factorization system $\langle L,R \rangle$ lifts to a factorization
system $\langle L',R' \rangle$ on $I'$, with $f \in R'$ iff $\pi(f)
\in R$ and $f$ is cartesian over $I$, and $f \in L'$ iff $\pi(f) \in
L$. In this case, we have $\Ar^{R'}(I') \cong \pi^*\Ar^R(I)$.
\end{exa}

While Example~\ref{cyl.exa} is the protopytical example of a
cofibration, it is Example~\ref{ar.exa}, Example~\ref{sq.exa} and
especially Example~\ref{facto.exa} that are very useful for actually
constructing fibrations and cofibrations. For another useful general
construction, consider a functor $\gamma:I' \to I$ between small
categories, and its left comma-fibers $I'/^\gamma i$, $i \in I$ of
Remark~\ref{comma.rem}. Then for any object $i \in I$, we have the
forgetful functor $p_i:I'/^\gamma i \to I'$, $\langle c,\alpha
\rangle \mapsto c$, and a map $f:i \to i'$ induces a functor
$f_!:I'/^\gamma i \to I'/^\gamma i'$, $\langle c,\alpha \rangle
\mapsto \langle c,f \circ \alpha \rangle$. For any fibration $\C' \to
I'$, we can define a fibration $\gamma_*\C' \to I$ with fibers
\begin{equation}\label{pd.fib}
(\gamma_*\C')_i = \Sec^\natural(I'/^\gamma i,p_i^*\C')
\end{equation}
and transition functors $(f_!)^*$. We then have a tautological
evaluation functor $\gamma^*\gamma_*\C' \to \C'$, cartesian over
$I'$, and for any fibration $\C \to I$, we have a tautological
functor $\C \to \gamma_*\gamma^*\C$, cartesian over $I$. This
defines a sort of a $2$-categorical adjunction between $\gamma^*$
and $\gamma_*$, in that for any fibrations $\C' \to I'$, $\C \to I$,
we have a natural equivalence of categories
$$
\Fun_{I'}^\natural(\gamma^*\C,\C') \cong
\Fun_I^\natural(\C,\gamma_*\C').
$$
If $\gamma$ is a cofibration, one can use the equivalences
\eqref{sec.I} to replace the left comma-fibers $I'/^\gamma i$ in
\eqref{pd.fib} with the usual fibers $I'_i$. In fact, in this case,
one can define a fibration $\gamma_{**}\C'$ with fibers
$(\gamma_{**}\C')_i \cong \Sec(I'_i,\C'_i)$, where $\C'_i$ is the
restriction of the fibration $\C'$ to the fiber $I'_i \subset I$, and
transition functors
\begin{equation}\label{st.st}
\begin{CD}
\Sec(I'_{i'},\C'_{i'}) @>{(f_!)^*}>> \Sec(I'_i,(f_!)^*\C'_{i'}) @>>>
\Sec(I'_i,\C'_i)
\end{CD}
\end{equation}
for any map $f:i \to i'$, where $f_!:I'_i \to I'_{i'}$ is the
transition functor of the cofibration $\gamma$, and the second
functor is induced by the transition functors of the fibration $\C'
\to I'$ along cocartesian liftings of the map $f$. Then $\gamma_*\C'
\subset \gamma_{**}\C'$ is the full subcategory spanned by
$\Sec^\natural(I'_i,\C'_i) \subset \Sec(I'_i,\C'_i)$.

\begin{remark}
We note that the $2$-functor $I^o \to \cCat$ corresponding to a
fibration $\gamma_*\C'$ of the form \eqref{pd.fib} is actually
strict. This can be used for rigidifying arbitrary pseudofunctors:
for any fibration $\C \to I$, we have $\C \cong \Id_*\C$, and
$\Id_*\C$ corresponds to a strict $2$-functor $I^o \to \cCat$ (in
particular, we have a genuine functor $I^o \to \Cat$). In effect,
one can consider the universal fibration $\Cat^\hdot \to \Cat^o$
corresponding to the tautological pseudofunctor $\Cat \to
\cCat$. Its fiber $\Cat^\hdot_{\C}$ over some $\C \in \Cat$ is $\C$
itself, and for any fibration $\C \to I$, there exists a functor
$\gamma:I \to \Cat^o$ and an equivalence $\eps:\C \cong
\gamma^*\Cat^\hdot$, and the pair $\langle \gamma,\eps \rangle$ is
unique up to a unique equivalence. However, these are equivalences
and not isomorphisms --- there is no way to state the Grothendieck
construction purely in terms of the category $\Cat$ without
introducing the notion of an equivalence. Replacing $\Cat$ with its
quotient $\bCat$ of Lemma~\ref{bcat.le} does not help either, since
a functor $I \to \bCat$ contains strictly less information that a
pseudofunctor $I \to \cCat$ (even if the latter is considered ``up to
an equivalence'').
\end{remark}

\begin{exa}\label{sch.exa}
Fix a Noetherian commutative ring $k$, denote by $\Comm(k)$ the
category of finitely generated unital associative commutative
$k$-algebras, and let $\Comm(k)\amod$ be the category of pairs
$\langle A,M \rangle$, $A \in \Comm(k)$, $M \in A\amod$ an
$A$-module, with maps $\langle A,M \rangle \to \langle A',M'
\rangle$ given by pairs $\langle f,g \rangle$ of an algebra map $f:A
\to A'$ and an abelian group map $g:M \to M'$ such that
$g(am)=f(a)g(m)$, $a \in A$, $m \in M$. Then the obvious forgetful
functor $\Comm(k)\amod \to \Comm(k)$ is bifibration, with fibers
$\Comm(k)\amod_A \cong A\amod$, cofibration transition functors $M
\mapsto A' \otimes_A M$ and the tautological fibration transition
functors sending an $A'$-module $M$ to itself with the induced
$A$-module structure.

More generally, let $\Sch(k)$ be the category of separated schemes
of finite type over $\Spec k$, let $j:\Comm(k)^o \to \Sch(k)$ be the
embedding $A \mapsto \Spec A$, and denote $\QCoh(\Sch(k)) =
(j_*\Comm(k)\amod^o)^o$. Then the cofibration $\QCoh(\Sch(k)) \to
\Sch(k)^o$ is again a bifibration, and for any $X \in \Sch(k)$, the
fiber $\QCoh(\Sch(k))_X \cong \QCoh(X)$ is the category of
quasicoherent sheaves on $X$. Explicitly, objects in
$\QCoh(\Sch(k))$ are pairs $\langle X,\F \rangle$, $X \in \Sch(k)$,
$\F \in \QCoh(X)$, and maps $\langle X',\F' \rangle \to \langle X,\F
\rangle$ are pairs of a scheme map $f:X \to X'$ and a sheaf map
$f^*\F \to \F'$ (note the change of variance). A scheme map $f$
defines a map $f^o$ in $\Sch(k)^o$, and the transition functors for
the bifibration $\QCoh(\Sch(k)) \to \Sch(k)^o$ are $f^{o*} \cong
f_*$ and $f^o_!  \cong f^*$.

Even more generally, we can consider complexes of quasicoherent
sheaves instead of sheaves themselves, denote by
$\QCoh_\idot(\Sch(k)) \to \Sch(k)^o$ the corresponding bifibration,
and restrict our attention to the full subcategory
$\QCoh_\idot^\flat(\Sch(k)) \subset \QCoh_\idot(\Sch(k))$ of
$h$-flat complexes (that is, complexes $\F_\idot$ such that
$\F_\idot \otimes_{\calo_X} -$ preserves quasiisomorphisms). Since
these are preserved by pullback functors $f^*$, the projection
$\QCoh_\idot^\flat(\Sch(k)) \to \Sch(k)^o$ is also a cofibration,
and since pullbacks preserve quasiisomorphisms of $h$-flat
complexes, we can then invert quasiisomorphisms in each fiber and
obtain a cofibration $\D(\Sch(k)) \to \Sch(k)^o$ with fibers
$\D(\Sch(k))_X \cong \D(X)$, the derived categories of
$\QCoh(X)$. Its transition functors are $f^o_! \cong f^*$, and the
isomorphisms \eqref{fib.eq} are induced by the corresponding
isomorphisms for the cofibration $\QCoh_\idot^\flat(\Sch(k)) \to
\Sch(k)^o$. This ``derived'' cofibration $\D(\Sch(k)) \to \Sch(k)^o$
is then again a bifibration, with fibration transition functors
given by the derived functors $R^\hdot f_*$. Alternatively, one can
start with $h$-injective complexes of flasque quasicoherent sheaves,
and functors $f^{o*} \cong f_*$; the result is the same.
\end{exa}

\section{Definitions and statements.}\label{2cat.sec}

\subsection{Nerves and $2$-categories.}

In practical applications, even if one considers $2$-functors and
not strict $2$-functors, the notion of a strict $2$-category is
still too strict. One wants to relax it by allowing composition
functors that are unital and associative only up to a fixed
isomorphism, subject to higher constraints, and then one has to
modify the notion of a $2$-functor to take account of these extra
pieces of data. A quick look at Definition~\ref{st.2fun.def} will
show that this is maybe not the best of ideas (at least if one is
intent on writing down complete proofs). We will now describe an
alternative formalism based of the notion of a ``nerve''.

As usual, we denote by $\Delta$ the category of ordinals $[n] =
\{0,\dots,n\}$, $n \geq 0$, and order-preserving maps between
them. Since ordinals can be thought of as small categories, we have
a fully faithful embedding $\delta:\Delta \subset \Cat$. We have
$[n]^o \cong [n]$, and this defines an involution $\iota:\Delta \to
\Delta$.  For any $n \geq l \geq 0$, we denote by $s,t:[l] \to [n]$
the embeddings identifying $[l]$ with an initial resp.\ terminal
segment of the ordinal $[n]$, and we note that we have a cocartesian
square
\begin{equation}\label{seg.sq}
\begin{CD}
[0] @>{t}>> [l]\\
@V{s}VV @VV{s}V\\
[n-l] @>{t}>> [n].
\end{CD}
\end{equation}
Recall that a {\em simplicial set\/} is a functor $X:\Delta^o \to
\Sets$, and the {\em nerve} $N(I)$ of a small category $I$ is a
simplicial set sending $[n]$ to the set of functors $[n] \to I$. In
particular, $0$-simplicies $i \in N(I)([0])$ are objects $i \in I$,
and $1$-simplices $f \in N(I)([1])$ are arrows $f:i \to i'$. The two
projections $N(I)(s),N(I)(t):N(I)([1]) \to N(I)([0])$ send an arrow
to its source resp.\ its target. A small category is completely
defined by its nerve, and a simplicial set $X$ is the nerve of a
small category iff it sends cocartesian squares \eqref{seg.sq} to
cartesian squares of sets (this is known as the {\em Segal
  condition}). One can also apply the Grothendieck construction to
the nerve; this gives a fibration $\Delta I \to \Delta$ with
discrete fibers $(\Delta I)_{[n]} = N(I)([n])$. The category $\Delta
I$ is called the {\em simplicial replacement} of the category $I$.

\begin{exa}\label{d.n.exa}
If one considers $[n] \in \Delta$ as a small category, then $\Delta
[n]$ is the left comma-fiber $\Delta /^{\Id}[n]$ of
Remark~\ref{comma.rem} of the identity functor --- that is, the
category of objects $[m] \in \Delta$ equipped with a map $[m] \to
[n]$. In particular, we have the tautological object $\id \in
(\Delta [n])_{[n]} \subset \Delta [n]$ that corresponds to the
identity map $\id:[n] \to [n]$.
\end{exa}

We now observe that the correspondence $I \mapsto \Delta I$ can be
immediately generalized to strict $2$-categories.

\begin{defn}\label{2repl.def}
The {\em simplicial replacement} $\Delta \C$ of a strict
$2$-category $\C$ is the fibration $\Delta \C \to \Delta$ with
fibers $(\Delta \C)_{[n]} = \Fun^2([n],\C)$, and transition functors
given by pullbacks with respect to maps $f:[n] \to [m]$.
\end{defn}

By Definition~\ref{fun2.def}, if $\C = I$ is a usual category, then
$\Fun^2([n],I)$ is exactly the set $N(I)([n])$, so that in this
case, Definition~\ref{2repl.def} reduces to our previous notion of a
simplicial replacement. We also note that any $2$-functor $\gamma:\C
\to \C'$ induces a functor $\Delta(\gamma):\Delta\C \to \Delta\C'$,
cartesian over $\Delta$, and as soon as $\C$ is small, we also have
a functor
\begin{equation}\label{fun2.eq}
\Fun^2(\C,\C') \to \Fun^\natural_\Delta(\Delta\C,\Delta\C').
\end{equation}
For any $\C$, the fiber $(\Delta\C)_{[0]}$ is discrete and consists
of objects of $\C$, while $(\Delta\C)_{[1]}$ is the disjoint union
of all the categories of morphisms $\C(-,-)$. Moreover, say that a
fibration $\E \to \Delta$ {\em satisfies the Segal condition} if for
any square \eqref{seg.sq}, the corresponding functor
\begin{equation}\label{seg.eq}
s^* \times t^*:\E_{[n]} \to \E_{[l]} \times_{\E_{[0]}} \E_{[n-l]}
\end{equation}
is an equivalence of categories. Then $\Delta\C$ always satisfies the
Segal condition. In particular, for $n=2$, $l=1$, the category
$(\Delta\C)_{[2]} \cong (\Delta\C)_{[1]} \times_{(\Delta\C)_{[0]}}
(\Delta\C)_{[1]}$ is the category of composable pairs of morphisms
in $\C$. The composition is then given by the transition functor
\begin{equation}\label{m.eq}
m^*:(\Delta\C)_{[2]} \cong (\Delta\C)_{[1]} \times_{(\Delta\C)_{[0]}}
(\Delta\C)_{[1]} \to (\Delta\C)_{[1]}
\end{equation}
corresponding to the embedding $m:[1] \to [2]$, $0 \mapsto 0$, $1
\mapsto 2$, and the identity morphisms are given by the embedding
\begin{equation}\label{e.eq}
e^*:(\Delta\C)_{[0]} \to (\Delta\C)_{[1]}
\end{equation}
corresponding to the unique projection $e:[1] \to [0]$.

\begin{lemma}\label{repl.le}
For any strict $2$-categories $\C$, $\C'$ with $\C$ small, the
functor \eqref{fun2.eq} is an equivalence.
\end{lemma}

\proof{} To obtain an inverse equivalence, note that a functor
$\gamma:\Delta\C \to \Delta\C'$ over $\Delta$ induces a map of
objects $\gamma([0])$, and a functor $\gamma([1])$ of the categories
of morphisms. If $\gamma$ is cartesian, then we also have the
isomorphisms \eqref{fib.fu} for the maps $m$, $e$ of \eqref{m.eq},
\eqref{e.eq}, and these provide the functorial isomorphisms of
Definition~\ref{st.2fun.def}. A map $\alpha:\gamma \to \gamma'$
between two cartesian functors has the component $\alpha([0])$ that
must be an identity map since $(\Delta\C')_{[0]}$ is discrete, and
$\alpha([1])$ gives a $2$-morphism in the sense of
Definition~\ref{fun2.def}.
\endproof

\begin{corr}\label{equi.corr}
A $2$-functor $\gamma$ is a strong equivalence in the sense of
Definition~\ref{stro.def} if and only if $\Delta(f)$ is an
equivalence of categories.
\end{corr}

\proof{} Clear. \endproof

\begin{defn}\label{2cat.def}
A {\em $2$-category} is a fibration $\C \to \Delta$ satisfying the
Segal condition \eqref{seg.eq} such that the fiber $\C_{[0]}$ is
discrete. A {\em $2$-functor} between $2$-categories $\C$, $\C'$ is
a functor $\gamma:\C \to \C'$ cartesian over $\Delta$.
\end{defn}

Motivated by Lemma~\ref{repl.le}, we denote $\Fun^2(\C,\C') =
\Fun^\natural_\Delta(\C,\C')$ for any $2$-categories $\C$, $\C'$ in
the sense of Definition~\ref{2cat.def}, and we let $\Fun^2(I,\C) =
\Fun^2(\Delta I,\C)$ for a small category $I$. We note that by virtue
of Example~\ref{d.n.exa}, we have a tautological equivalence
$\Y:\Fun^2([n],\C) \cong \C_{[n]}$ sending a $2$-functor $\gamma$ to
its value on $\id \in (\Delta [n])_{[n]}$. For any strict
$2$-category $\C$, the simplicial replacement $\Delta\C$ is a
$2$-category in the sense of Definition~\ref{2cat.def}, and for any
$2$-functor $\gamma:\C \to \C'$, the functor
$\Delta(\gamma):\Delta\C \to \Delta\C'$ is a $2$-functor. An
arbitrary $2$-category $\C$ in the sense of
Definition~\ref{2cat.def} still has objects $c \in \C_{[0]}$, and
just as in the strict case, the projection $s^* \times t^*:\C_{[1]}
\to \C_{[0]}$ provides a decompositon
\begin{equation}\label{c.1.eq}
\C_{[1]} = \coprod_{c,c' \in \C_{[0]}}\C(c,c')
\end{equation}
into categories of morphisms. Then \eqref{m.eq} and \eqref{e.eq}
provide composition functors and identity objects. However, in
general, all the axioms only hold up to canonical isomorphisms
induced by the maps \eqref{fib.fu} for the structural fibration $\C
\to \Delta$.

For any $2$-category $\C$, the $2$-opposite category is
$\C^\tau=\C^o_\perp$, and the $1$-opposite $2$-category is $\C^\iota
= \iota^*\C$, where $\iota:\Delta \to \Delta$ is the involution $[n]
\mapsto [n]^o$. The {\em $2$-product} of $2$-categories $\C$, $\C'$
is given by $\C \times^2 \C' = \C \times_{\Delta} \C'$. Note that
any object $c \in \C'_{[0]}$ uniquely extends to a cartesian section
$\iota(c):\Delta \to \C'$ of the fibration $\C' \to \Delta$, and
this gives a full embedding $\id \times \iota(c):\C \to \C \times^2
\C'$ whose essential image is denoted by $\C \times \{c\} \subset
\C'$.

\begin{remark}\label{cova.rem}
Definition~\ref{2cat.def} is essentially the same as
\cite[Definition 6.1]{ka.hh} but there is one difference:
\cite{ka.hh} describes a $2$-category $\C$ of
Definition~\ref{2cat.def} in terms of the transpose cofibration
$\C_\perp \to \Delta^o$ rather then the fibration $\C \to
\Delta$. The two approaches are equivalent but one has to choose
one. It seems that fibrations make for simpler formulas when dealing
with sets and geometric objects, while cofibrations work better for
rings, algebras and the like. Since \cite{ka.hh} deals with algebra
rather than geometry, cofibrations were chosen there, and this is
the choice that would have been better for Section~\ref{app.sec}
below. However, for the rest of this Section and for
Section~\ref{adj.sec} fibrations are more convenient, so this is
what we use.
\end{remark}

\subsection{Grothendieck construction for
  $2$-categories.}\label{2groth.subs}

As a first example for Definition~\ref{2repl.def}, let us describe
the simplicial replacement $\Delta\cCat$ of the $2$-category
$\cCat$. We first observe that any fibration $\E \to \Delta$ that
satisfies the Segal condition can be turned into a $2$-category in
the following way. Take any fibration $\E \to \Delta$, let
$\eps:\ppt \to \Delta$ be the embedding onto $[0] \in \Delta$, let
$\E_0 = \E_{[0]} = \eps^*\E$, and let $\E'_0 = \E_{0,\Id} \subset
\E_0$ be its underlying discrete subcategory. Define the {\em
  reduction} $\E^{red}$ of the fibration $\E$ by the cartesian
square
\begin{equation}\label{red.sq}
\begin{CD}
\E^{red} @>>> \E\\
@VVV @VVV\\
\eps_*\E'_0 @>>> \eps_*\E_0.
\end{CD}
\end{equation}
Then the square is fibered over $\Delta$, and by \eqref{pd.fib}, we
have $(\eps_*\E_0)_{[n]} \cong \E_0^{n+1}$ for any $[n] \in \Delta$,
and similarly for $\E'_0$. Thus both categories in the bottom line
of \eqref{red.sq} satisfy the Segal condition. Then if so does $\E$,
the same holds for $\E^{red}$, and since $\E^{red}_{[0]} = \E'_0$ is
discrete, it is a $2$-category in the sense of
Definition~\ref{2cat.def}.

Now let $c$ be the class of all cofibrations in $\Cat$, consider the
fibration $t:\Ar^c(\Cat) \to \Cat$ of Example~\ref{sq.exa} (that
applies by Remark~\ref{iso.rem}), and let $\E = \delta^*\Ar^c(\Cat)$
be its restriction with respect to the full embedding $\delta:\Delta
\to \Cat$. We claim that
\begin{equation}\label{ccat.eq}
\Delta\cCat \cong \E^{red}.
\end{equation}
Indeed, by definition, for any $[n] \in \Delta$, $\E_{[n]} \cong
\Cat_\natural/[n]$ is the category of small categories cofibered
over $[n]$, and strict cocartesian functors between them. By the
Grothendieck construction, these categories correspond to
$2$-functors $[n] \to \cCat$. Morphisms between some $\C,\C' \in
\E_{[n]}$ are strict cocartesian functors $\gamma:\C \to \C'$ over
$[n]$, but morphisms in the reduction $\E^{red}_{[n]}$ are functors
such that $\gamma(l)=\id$ for any $l \in [n]$. The only remaining
data are then the morphisms \eqref{fib.fu} -- or rather, their
cofibration analogs -- and these correspond on the nose to
$2$-morphisms of Definition~\ref{fun2.def}.

\begin{remark}
The reader might wonder what changes if instead of strict
cocartesian functors we consider all cocartesian functors. The
answer is ``in the end, nothing'': while $\E$ would change, its
reduction would be the same.
\end{remark}

As an application of the identification \eqref{ccat.eq}, let us
describe a version the Grothendieck construction of
Subsection~\ref{groth.subs} with $I$ replaced by an arbitrary
$2$-category $\C$. Say that a map $f:[m] \to [n]$ in $\Delta$ is
{\em special} if $f(0)=0$. Denote by $+$ the class of all special
maps, with $\Delta_+ \subset \Delta$ being the corresponding dense
subcategory, and let $t$ be the class of all terminal embeddings
$t:[m] \to [n]$. Then $\langle +,t \rangle$ is a factorization
system on $\Delta$, and for any $2$-category $\C$, we have its
lifting $\langle +,t \rangle$ of Example~\ref{facto.exa}. Say that a
map $f$ in $\C$ is {\em special} if $f \in +$ (that is, $f$ goes to
a special map under the fibration $\C \to \Delta$).

\begin{defn}
A fibration $\C' \to \C$ over a $2$-category $\C$ is {\em special}
if for any special map $f$ in $\C$, the transition functor $f^*$ is
an equivalence.
\end{defn}

In the prototypical example of a special fibration, $\C$ is
$\Delta\cCat$, and the fibration $\Delta^\hdot\cCat \to \Delta\cCat$
has fibers
\begin{equation}\label{del.cat}
\Delta^\hdot\cCat_{\C} = \Sec_\natural([n],\C),
\end{equation}
where the cofibration $\C \to [n]$ corresponds to an object in
$\Delta\cCat$ via \eqref{ccat.eq}. The transition functors are given
by pullbacks. Since $0 \in [n]$ is the initial object, we have the
equivalence \eqref{sec.I} -- or rather, its analog for cofibrations
and cocartesian sections -- and in our case, it reads as
$\Sec_\natural([n],\C) \cong \C_0$. Therefore the fibration
\eqref{del.cat} is special. It is also universal, in the following
sense.

\begin{lemma}\label{2groth.le}
For any special fibration $\C' \to \C$ over a $2$-category $\C$,
there exists a $2$-functor $\gamma:\C \to \Delta\cCat$ such that
$\C' \cong \gamma^*\Delta^\hdot\cCat$.
\end{lemma}

\proof{} Note that for any $[n] \in \Delta$, the fiber
$\Ar^t(\Delta)_{[n]}$ of the cofibration $t$ of
Example~\ref{facto.exa} is naturally identified with $[n]^o$, by
sending $t:[m] \to [n]$ to $t(0) \in [n]$, and the same is then true
for the fibers of the cofibration $t:\Ar^t(\C) \to \C$. For any $c
\in \C$, the projection $s:\Ar^t(\C) \to \C$ restricts to a functor
$s_c:[n]^o \cong \Ar^t(\C)_c \to \C$, so we have the fibration
$s_c^*\C' \to [n]^o$. The functor $\gamma$ sends $c$ to the
transpose cofibration $(s_c^*\C')_\perp \to [n]$.
\endproof

\begin{remark}\label{del.dot.rem}
The fibration $t^o:\Ar^t(\Delta)^o \to \Delta^o$ opposite to the
cofibration $t$ corresponds to the tautological embedding $\Delta
\subset \Cat$ by the Grothendieck construction. Another
interpretation is that $\Ar^t(\Delta)$ is the category
$\Delta^\hdot$ of pairs $\langle [n],l\rangle$, $[n] \in \Delta$, $l
\in [n]$, with maps from $\langle [n],l \rangle$ to $\langle [n'],l'
\rangle$ given by maps $f:[n] \to [n']$ such that $f(l) \geq
l'$. Then $t$ is the forgetful functor $\langle [n],l \rangle
\mapsto [n]$, and $s:\Ar^t(\Delta) \to \Delta$ sends $\langle [n],l
\rangle$ to $[n-l]$. This interpretation is quite useful in many
constructions involving nerves. For example, the simplicial
replacement $\Delta I$ of a category $I$ is given by
$$
\Delta I = (t_{**}(I \times \Ar^t(\Delta)))^{red},
$$
where $t_{**}$ is as in \eqref{st.st}, and $I \times \Ar^t(\Delta)
\to \Ar^t(\Delta)$ is the constant fibration with fiber $I$.
\end{remark}

\begin{remark}\label{red.rem}
Ideally, we would like to say that $\gamma$ in Lemma~\ref{2groth.le}
is unique up to an isomorphism, but this is not true. The reason for
this is the reduction procedure used in \eqref{ccat.eq}: by
definition, objects in $(\Delta\cCat)_{[0]}$ are small categories on
the nose, and two isomorphic but different categories stop being
isomorphic as objects in $\Delta\cCat$. The real way to get rid of
artefacts of this type would be to relax the condition of
Definition~\ref{2cat.def} and allow $\C_{[0]}$ to be a groupoid;
however, this would take one rather far away from the standard
theory of $2$-categories. As a temporary fix, one can choose a
specific model for $\Cat$ and similar categories where there is only
one object in each isomorphism class. Then $\gamma$ in
Lemma~\ref{2groth.le} becomes unique up to a unique isomorphism.
\end{remark}

If we fix a small category $I$, and fix a concrete model for $\Cat$
as in Remark~\ref{red.rem}, that the equivalence \eqref{fi.co.eq}
extends to a strong equivalence
$$
\cCat^\natural/I \cong \cCat_\natural/I^o,
$$ and all of the material in this Subsection immediately
generalizes to these strongly equivalent strict $2$-categories, with
the same proofs. Namely, for any $2$-category $\C$, a fibration $\C'
\to I \times \C$ is {\em special} if so is its composition with the
projection $I \times \C \to \C$. Objects in the simplicial
replacement $\Delta(\cCat_\natural/I^o)$ are pairs $\langle [n],\C
\rangle$ of an ordinal $[n] \in \Delta$ and a small cofibration $\C
\to I^o \times [n]$, and one defines a special fibration
$\Delta^\hdot(\cCat_\natural/I^o) \to I \times
\Delta(\cCat_\natural/I^o)$ by setting
$\Delta^\hdot(\cCat_\natural/I^o)_\C = \pi_*\C^\perp$, where
$\C^\perp \to I \times [n]^o$ is the transpose fibration to $\C$,
and $\pi:I \times [n]^o \to I$ is the projection. Then exactly the
same argument as in Lemma~\ref{2groth.le} shows that for any
$2$-category $\C$ equipped with a special fibration $\C' \to I
\times \C$, there exists a $2$-functor $\gamma:\C \to
\Delta(\cCat_\natural/I^o)$ and an equivalence $\C' \cong (\gamma
\times \Id)^*\Delta^\hdot(\cCat_\natural/I^o)$, with the uniqueness
properties described in Remark~\ref{red.rem}.

\subsection{Rigidification.}

Let us now show that our Definition~\ref{2cat.def} is equivalent to
the usual definition of weak $2$-categories and $2$-functors that
uses associativity isomorphisms and suchlike. Fortunately, it has
already been proved in \cite{ben} that any weak $2$-category in the
usual sense is equivalent to a strict one, so it suffices to
identify $2$-categories and $2$-functors of
Definition~\ref{2cat.def} with strict $2$-categories and
$2$-functors. For $2$-functors, this is Lemma~\ref{repl.le}, and the
result for $2$-categories is as follows.

\begin{theorem}\label{rigi.thm}
For any small $2$-category $\C$ in the sense of
Definition~\ref{2cat.def}, there exists a small strict $2$-category
$\R(\C)$ and an equivalence $\C \cong \Delta\R(\C)$, and $\R(\C)$ is
unique up to a strong equivalence in the sense of
Definition~\ref{stro.def}.
\end{theorem}

The idea of the proof of Theorem~\ref{rigi.thm} is to use the Yoneda
embedding. We consider the $1$-opposite $2$-category $\C^\iota$, and
to any object $c \in \C_{[0]}$, we associate a representable
$2$-functor $\C^\iota \to \cCat$. We then apply the Grothendieck
construction of Lemma~\ref{2groth.le} to encode these $2$-functors
by special fibrations $\C^\iota(c) \to \C^\iota$, we get a full
embedding $\C \to \Delta(\cCat^\natural/\C^\iota)$, and
$\cCat^\natural/\C^\iota$ is strict. In order to make this work, we
first need a $2$-categorical version of the Yoneda Lemma.

Let $\rho:\Delta_+ \to \Delta$ be the embedding functor, and note
that it has a left-adjoint $\lambda:\Delta \to \Delta_+$ that adds a
new initial element to an ordinal $[n]$. The composition $\kappa =
\rho \circ \lambda:\Delta \to \Delta$ sends an ordinal $[n]$ to
$[n]^< \cong [n+1]$ and comes equipped with the adjunction map
$a:\Id \to \kappa$ (given by the terminal segment embeddings $t:[n]
\to [n+1]$). By Remark~\ref{comma.rem}, since $[0] \in \Delta_+$ is
the initial object, any fibration $\E \to \Delta_+$ comes equipped
with a functor $\E \to \E_{[0]}$ right-adjoint to the embedding
$\E_{[0]} \to \E$. In particular, for any $2$-category $\C$, we can
take $\C_+ = \rho^*\C \to \Delta_+$, and we obtain a decomposition
\begin{equation}\label{C.facto}
\C_+ \cong \coprod_{c \in \C_{[0]}}\C_+(c).
\end{equation}
Applying $\lambda^*$, we obtain a decomposition of the
fibration $\kappa^*\C \cong \lambda^*\C_+$ into components $\C(c) =
\lambda^*\C_+(c)$, $c \in \C_{[0]}$, and for any $c$, the transition
functor $a^*$ provides a projection
\begin{equation}\label{C.c.eq}
\pi = a^*:\C(c) \to \C.
\end{equation}
The Segal condition immediately shows that this is a
fibration. Moreover, any object $c' \in \C_{[n]} \subset \C$ defines
an object $s^*c' \in \C_{[0]}$, where $s:[0] \to [n]$ is the initial
embedding, and again by the Segal condition, the fibers of the
fibration \eqref{C.c.eq} are given by
\begin{equation}\label{s.c}
\C(c)_{c'} \cong \C(c,s^*c'),
\end{equation}
where $\C(-,-)$ are the components of the decomposition
\eqref{c.1.eq}. In particular, the fibration \eqref{C.c.eq} is
special. Moreover, we have a distinguished object $\id_c \in \C(c,c)
\cong \C(c)_c \subset \C(c)$.

\begin{lemma}\label{yo.le}
For any small $2$-category $\C$ and $c,c' \in \C_{[0]}$, the functor
\begin{equation}\label{yo.eq}
\Y:\Fun^\natural_\C(\C(c),\C(c')) \to \C(c',c), \qquad \gamma
\mapsto \gamma(\id_c)
\end{equation}
is an equivalence of categories.
\end{lemma}

\proof{} The embedding $\lambda_c:\C(c) \cong \lambda^*\C_+(c) \to
\C_+(c)$ admits a right-adjoint functor $\rho_c:\C_+(c) \to \C(c)$
whose composition $\pi_c = \pi \circ \rho_c:\C_+(c) \to \C$ with the
fibration $\pi$ of \eqref{C.c.eq} is obtained by restricting the
natural embedding $\rho:\C_+ = \rho^*\C \to \C$ to the corresponding
component of \eqref{C.facto}. For any map $f$ in $\C_+(c)$, the map
$\rho_c(f)$ in $\C(c)$ is cartesian with respect to $\pi$, and for
any object $x \in \C(c)$, so is the adjunction map $a_c(x):x \to
\kappa_c(x)$, where $\kappa_c = \rho_c \circ \lambda_c$. Therefore
for any functor $F:\C(c) \to \C(c')$ cartesian over $\C$, and any
object $x \in \C(c)$, we have $F(x) \cong
\pi(a_c(x))^*F(\kappa_c(x))$, and the isomorphism is functorial with
respect to $F$ and $x$. This means that we have an equivalence
\begin{equation}\label{yo.yo.eq}
\Fun^\natural_\C(\C(c),\C(c')) \cong
\Sec^\natural(\C_+(c),\pi_c^*\C(c')), \qquad F \mapsto
\rho_c^*F,
\end{equation}
and it remains to observe that since $\pi_c:\C_+(c) \to \C$ sends
all maps to special maps, and the fibration $\C(c') \to \C$ is
special, all the transition functors of the fibration
$\pi_c^*\C(c')$ are equivalences. Therefore it is a
bifibration, its cartesian sections coincide with cocartesian ones,
and by the cofibration verison of \eqref{sec.I}, the right-hand side of
\eqref{yo.yo.eq} is equivalent to the fiber $\C(c',c) \cong
\pi_c^*\C(c')_o$ over the initial object $o \in \C_+(c)$.
\endproof

\proof[Proof of Theorem~\ref{rigi.thm}.] Uniqueness immediately
follows from Corollary~\ref{equi.corr}. To prove existence, say that
the {\em concatenation} $[l] \circ [n]$ of any two ordinals $[l],[n]
\in \Delta$ is their disjoint union $[l] \copr [n]$ ordered
left-to-right, and denote by $m:\Delta^2 \to \Delta$ the functor
sending $[l] \times [n]$ to $[l] \circ [n]$. Moreover, let
$p_0,p_1:\Delta^2 \to \Delta$ be the projections onto the first and
second factor, and note that we have natural maps $a_0:p_0 \to m$,
$a_1:p_1 \to m$. Let $\C^{(2)} = m^*\C$, and note that the
transition functors $a_0^*$, $a_1^*$ induce projections
$\pi_0,\pi_1:\C^{(2)} \to \C$. Consider the product
\begin{equation}\label{C.2.eq}
\pi = (\iota \circ \pi_0) \times \pi_1:\C^{(2)} \to \C^\iota \times
\C.
\end{equation}
Then the Segal condition immediately shows that this is a
fibration. Moreover, its composition $\C^{(2)} \to \C$ with the
projection $\C^\iota \times \C \to \C$ is also a fibration with
fibers $\C^{(2)}_c \cong \C^\iota(s^*c)$, where $s^*c$ is as in
\eqref{s.c}. Therefore the fibration \eqref{C.2.eq} is special, and
if we denote by $\R(\C) \subset \cCat^\natural/\C^\iota$ the full
$2$-subcategory spanned by special fibrations $\C^\iota(c) \to
\C^\iota$, $c \in \C^\iota_{[0]} = \C_{[0]}$, then
Lemma~\ref{2groth.le} provides a $2$-functor
\begin{equation}\label{Y.eq}
\Y:\C \to \Delta\R(\C).
\end{equation}
By construction, this functor is an identity over $[0] \in \Delta$,
and for any objects $c,c' \in \C_{[0]}$, its restriction $\Y(c,c')$
to the component $\C(c,c')$ of the decomposition \eqref{c.1.eq} is
inverse to the equivalence \eqref{yo.eq} of
Lemma~\ref{yo.le}. Therefore $\Y(c,c')$ is also an equivalence, so
that \eqref{Y.eq} is an equivalence over $[1]$. By the Segal
condition, it is an equivalence.
\endproof

\begin{remark}\label{yo.rem}
The main reason we restrict our attention to small categories in
Lemma~\ref{yo.le} and Theorem~\ref{rigi.thm} is to insure that
everything is well-defined --- {\em a priori}, there could be more
than a set of morphisms between two functors between large
categories. However, it is easy to see that even for a large $\C$,
Lemma~\ref{yo.le} works with the same proof, and {\em a posteriori},
functors from $\C(c)$ to $\C(c')$ cartesian over $\C$ form a
well-defined category. The Yoneda embedding of
Theorem~\ref{rigi.thm} works for large $2$-categories, too.
\end{remark}

\subsection{Monoidal structures.}\label{mon.subs}

Let us now consider a special class of $2$-ca\-te\-gories --- those
that only have one object. These correspond to unital monoidal
structures.

\begin{defn}\label{mon.def}
A {\em unital monoidal structure} on a category $\C$ is given by a
$2$-category $B\C$ with $B\C_{[0]} \cong \ppt$, equipped with an
equivalence $B\C_{[1]} \cong \C$. A unital monoidal structure on a
functor $\gamma:\C \to \C'$ between categories $\C$, $\C'$ with
unital monoidal structures $B\C$, $B\C'$ is a $2$-functor
$B\gamma:B\C \to B\C'$ that restricts to $\gamma$ over $[1]$.
\end{defn}

In this definition, $2$-categories and $2$-functors are understood
in the sense of Definition~\ref{2cat.def}. One can also consider
{\em strict unital monoidal categories} by requiring $B\C$ to be a
strict $2$-category, but this is not very useful in practical
applications. Some very simple monoidal structures that exist in
nature are indeed strict; one such is the concatenation product used
in the proof of Theorem~\ref{rigi.thm} (it has to be strict by
necessity, since the only isomorphisms in $\Delta$ are the identity
maps). However, even this simple example admits a more efficient
description in terms of Definition~\ref{2cat.def}. To construct the
corresponding $2$-category $B\Delta$, we need some simplicial
combinatorics.

Recall from Subsection~\ref{2groth.subs} that we denote by $+$ the
class of special maps in $\Delta$, and that the embedding
$\rho:\Delta_+ \to \Delta$ admits a left-adjoint $\lambda:\Delta \to
\Delta_+$, $[n] \mapsto [n]^<$. Note that $\lambda$ extends to the
category $\Delta^<$, where we interpret the new initial object $o
\in \Delta^<$ as the empty ordinal. Dually, say that a map $f:[l]
\to [n]$ is {\em antispecial} if $\iota(f)$ is special (or
equivalently, $f(l)=n$), let $-$ be the class of antispecial maps,
and note that $\iota:\Delta \to \Delta$ induces an equivalence
$\Delta_+ \cong \Delta_-$, while the embedding $\rho_\iota = \iota
\circ \rho \circ \iota:\Delta_- \to \Delta$ admits a left-adjoint
$\lambda_\iota = \iota \circ \lambda \circ \iota:\Delta \to
\Delta_+$, $[n] \mapsto [n]^>$ that also extends to
$\Delta^<$. Moreover, say that a map $f$ is {\em bispecial} if it is
both special and antispecial, let $\pm = + \cap -$, and note that
the embedding $\rho_\flat = \rho \circ \rho_\iota:\Delta_\pm \to
\Delta$ admits a left-adjoint $\lambda_\flat:\Delta \to
\Delta_{\pm}$ that also extends to $\Delta^<$. Note that $[0] \in
\Delta$ is both the initial and the terminal object both for
$\Delta_+$ and $\Delta_-$ and the terminal object for $\Delta_\pm$,
while the initial object for $\Delta_\pm$ is $[1]$.

Next, recall that we have the factorization system $\langle +,t
\rangle$ on $\Delta$, where $t$ stands for the class of terminal
embeddings $t:[l] \to [n]$, and dually, we have the factorization
system $\langle -,s \rangle$, where $s$ is the class of the initial
embeddings $s:[l] \to [n]$. Say that a map $f$ in $\Delta$ is an
{\em anchor map} if it decomposes as $f = s \circ t$, with $s$
resp.\ $t$ an initial resp.\ a terminal embedding, or equivalently,
$f:[m] \to [n]$ is injective and identifies the ordinal $[m]$ with
some segment $\{l,l+1,\dots,l+m\} \subset [n]$ of the ordinal
$[n]$. Then we also have a factorization system $\langle \pm,a
\rangle$ on $\Delta$, where $a$ stands for the class of anchor maps.

Finally, observe that if we treat a map $f:[m] \to [n]$ in $\Delta$
as a functor between small categories, then a right-adjoint
$f^\dg:[n] \to [m]$ is by definition given by $f^\dg(l) = \max\{l'
\in [m]|f(l') \leq l\}$.  This is well-defined if and only if all
the sets in the right-hand side are non-empty, or equivalently,
$f(0) \leq 0$, so that $f$ must be special. The adjoint $f^\dg$ is
then automatically anti-special, and we obtain equivalences
\begin{equation}\label{p.m.o}
\Delta_+ \cong \Delta_-^o, \qquad \Delta_+^o \cong \Delta_-.
\end{equation}
Moreover, $f^\dg$ is special if and only if $f^{-1}(0) = 0$, and this
means that $f$ lies in the image of the functor $\lambda:\Delta^<
\to \Delta$. Therefore \eqref{p.m.o} induces equivalences
\begin{equation}\label{jo.eq}
\Delta_\pm \cong \Delta^{<o}, \qquad \Delta_\pm^o \cong \Delta^<.
\end{equation}
This observation is sometimes called the {\em Joyal duality},
although we have not been able to trace the origin of the name (nor
of the observation).

Now, the concatenation product $- \circ -$ obviosly extends to the
category $\Delta^<$, and it turns out that it is simpler to first
describe the corresponding monoidal structure on the opposite
category $\Delta^{<o}$. We have $B\Delta^{<o} \cong
\Ar^\pm(\Delta)$, with the fibration $s:\Ar^\pm(\Delta) \to \Delta$
of Example~\ref{facto.exa}, and the identification
$B\Delta^{<o}_{[1]} = \Delta_\pm \cong \Delta^{<o}$ is
\eqref{jo.eq}. In terms of the bispecial category $\Delta_\pm$, the
product is given by the {\em reduced concatenation} $[n] * [m] = [n]
\copr_{[0]} [m]$, where the coproduct is taken with respect to the
embeddings $t:[0] \to [n]$, $s:[0] \to [m]$ (so that what we
consider is a cocartesian square \eqref{seg.sq}).

Any monoidal structure $B\C$ on a category $\C$ induces a monoidal
structure $B\C^o = (B\C)^\tau = (B\C)_\perp^o$ on the opposite
category $\C^o$, so that we also obtain the monoidal structure on
$\Delta^<$. However, it also has a more direct description in terms
of the embedding $\lambda_\flat:\Delta^< \to \Delta_\pm$. Namely,
note that this embedding is injective on both objects and morphisms,
and its image consists of objects $[n] \in \Delta_\pm$ such that the
unique bispecial map $[1] \to [n]$ is injective (that is, $n \geq
1$), and bispecial maps $f:[n] \to [m]$ such that $f^{-1}(0)=0$ and
$f^{-1}(n) = m$. We can then consider the subcategory $\Ar^c(\Delta)
\subset \Ar(\Delta)$ spanned by injective bispecial arrows and maps
between them represented by commutative squares \eqref{ar.sq} that
are also cartesian. One checks that the projection $s:\Ar^c(\Delta)
\to \Delta$ is a fibration, and we have $B\Delta^< \cong
\Ar^c(\Delta)$.

Since concatenation and reduced concatenation are different, the
functors $\rho_\flat:\Delta_\pm \to \Delta \subset \Delta^<$,
$\lambda_\flat:\Delta^< \to \Delta_\pm$ are not monoidal. However,
we do have an obvious functorial map $[m] \circ [n] \to [m] * [n]$,
$[n],[m] \in \Delta$ from the usual concatenation to the reduced
one, and this is an example of the following useful $2$-categorical
structure.

\begin{defn}\label{lax.def}
For any $2$-categories $\C_0$, $\C_1$, a {\em co-lax $2$-functor} from
$\C_0$ to $\C_1$ is a functor $\gamma:\C_0 \to \C_1$ over $\Delta$ that
is cartesian over all anchor maps, and a {\em lax $2$-functor} is a
co-lax $2$-functor $\gamma:\C_0^\tau \to \C_1^\tau$ between the
$2$-opposite $2$-categories.
\end{defn}

For strict $2$-categories, one can describe lax and co-lax
$2$-functors in terms of Definition~\ref{st.2fun.def}: for a lax
$2$-functor $\gamma$, one allows the maps $\gamma_c$,
$\gamma_{c,c',c''}$ that are not invertible, and for a co-lax
$\gamma$, they also go in the opposite direction. In terms of
Definition~\ref{2cat.def}, these maps are the maps \eqref{fib.fu}
for the transition functors \eqref{e.eq} and \eqref{m.eq}. A co-lax
monoidal structure on a functor $\gamma:\C \to \C'$ between monoidal
categories $\C$, $\C'$ is a co-lax $2$-functor $B\gamma:B\C \to
B\C'$, and a lax monoidal structure on $\gamma$ is a co-lax one on
$\gamma^o$.

If $\gamma:\C_0 \to \C_1$ is a co-lax $2$-functor such that
$\gamma([0])$ is an isomorphism and $\gamma([1])$ admits a
right-adjoint $\gamma([1])^\dg$, then for any $[n] \in \Delta$,
$\gamma([n])$ admits a right-adjoint $\gamma([n])^\dg$ by virtue of
the Segal condition. In such a situation, Example~\ref{adj.perp.exa}
provides a functor $\gamma^\dg_\perp$ whose opposite
$\gamma^\dg:\C_1^\tau \to \C_0^\tau$ is a lax $2$-functor from
$\C_1$ to $\C_0$. In particular, if we have monoidal categories
$\C$, $\C'$ and functor $\gamma:\C \to \C'$ that admits a
right-adjoint $\gamma^\dg$, then a co-lax structure $B\gamma$ on
$\gamma$ induces a lax structure $B\gamma^\dg = (B\gamma)^\dg$ on
$\gamma^\dg$, and vice versa.

With these definitions, the embedding $\lambda_\flat:\Delta^< \to
\Delta_\pm$ has a natural co-lax monoidal structure given by the
embedding $\Ar^c(\Delta) \subset \Ar^b(\Delta)$, and the
right-adjoint embedding $\rho_\flat:\Delta_\pm \to \Delta \subset
\Delta^<$ then has a lax monoidal structure by adjunction. However,
the two structures are actually the same, since under the
identifications \eqref{jo.eq}, we have $\rho_\flat \cong
\lambda_\flat^o$.

\section{Adjunction.}\label{adj.sec}

\subsection{Adjoint pairs.}

For any $2$-category $\C$ and two objects $c,c' \in \C_{[0]}$, an
{\em adjoint pair} of maps between $c$ and $c'$ is a quadruple
$\langle f,f^\vee,a,a^\vee \rangle$ of objects $f \in \C(c,c')$,
$f^\vee \in \C(c',c)$ and maps $a:f^\vee \circ f \to \id_c$,
$a^\vee:\id_{c'} \to f \circ f^\vee$, subject to the relations
\begin{equation}\label{adj.rel}
(\id_f \circ a) \circ (a^\vee \circ \id_f) = \id_f, \quad (a \circ
  \id_{f^\vee}) \circ (\id_{f^\vee} \circ a^\vee) = \id_{f^\vee}.
\end{equation}
If $\C=\cCat$, then this is the usual definition of an adjoint pair
of functors. For any $\C$, we say that $f \in \C(c,c')$ is {\em
  reflexive} if it extends to an adjoint pair, and we say that a
reflexive morphism $f \in \C(c,c')$ is an {\em equivalence} if there
exists an adjoint pair $\langle f,f^\vee,a,a^\vee \rangle$ with
invertible $a$ and $a^\vee$. We note that all the usual properties
of adjunction in $\cCat$ extend to an arbitrary $\C$: reflexivity is
closed under compositions, the adjoint $f^\vee$ to a reflexive $f$
is unique up to a unique isomorphism, and $f$ is an equivalence if
and only if there exists $f^\vee \in \C(c',c)$ such that $\id_c
\cong f^\vee \circ f$ and $\id_{c'} \cong f \circ f^\vee$ (to reduce
the general case to the case $\C=\cCat$, one can use the Yoneda
embedding of Theorem~\ref{rigi.thm}, and this also works for large
$\C$ by Remark~\ref{yo.rem}).

\begin{remark}
The above definition is almost the same as in \cite[Subsection
  7.1]{ka.hh}, but with one difference: here $f^\vee$ is actually
left-adjoint to $f$, while in \cite{ka.hh} it is right-adjoint. The
reason for the discrepancy is the change of variance mentioned in
Remark~\ref{cova.rem}.
\end{remark}

By Example~\ref{d.n.exa}, we have $\Fun^2([1],\C) \cong \C_{[1]}$
for any $2$-category $\C$, so that morphisms in $\C$ correspond
bijectively to $2$-functors $\nat = \Delta [1] \to \C$. Analogously,
equivalences correspond to $2$-functors $\eq = \Delta e(\{0,1\}) \to
\C$. The tautological embedding $[1] \to e(\{0,1\})$ defines a
$2$-functor $\beta:\nat \to \eq$, and a morphism is an equivalence
if and only if the corresponding $2$-functor $\nat \to \C$ factors
through $\beta$. It turns out that one can reasonably explicitly
construct
a factorization
\begin{equation}\label{nat.adj.eq}
\begin{CD}
\nat @>{\delta}>> \adj @>{\nu}>> \eq
\end{CD}
\end{equation}
of the $2$-functor $\nat \to \eq$ and a $2$-category $\adj$ that
classifies adjoint pairs in the same way.

To do this, note that by definition, $\eq = \Delta e(\{0,1\})$ is
the category of pairs $\langle [n],e \rangle$, $[n] \in \Delta$,
$e:[n] \to e(\{0,1\})$ a functor (or equivalently, a map of sets
$[n] \to \{0,1\}$). Then we have functorial subsets $[n]_l =
e^{-1}(l) \subset [n]$, $l=0,1$. Say that a map $f:[n] \to [m]$ in
$\eq$ is {\em $l$-special}, $l=0,1$, if $f:[n]_l \to [m]_l$ is an
isomorphism. Then the classes $0$ and $1$ of $0$-special and
$1$-special maps in $\eq$ are closed and form a factorization system
in $\eq$ in the sense of Example~\ref{facto.exa}, in either order
(for a formal proof of this, see \cite[Lemma 7.4]{ka.hh}). Moreover,
if $f$ is bispecial, then so are both of the components of the
corresponding decomposition \eqref{facto.dia}. Then the
anchor/bispecial factorization system on $\Delta$ lifts to $\eq$ by
Example~\ref{facto.exa}, and any map $f:[n] \to [m]$ in $\eq$
uniquely factors as
\begin{equation}\label{3.eq}
\begin{CD}
[n] @>{b_1}>> [n'] @>{b_0}>> [n''] @>{a}>> [m],
\end{CD}
\end{equation}
with bispecial $1$-special $b_1$, bispecial $0$-special $b_0$, and
anchor $a$. This means that in particular, the class $\pm \cap 1$ of
bispecial $1$-special maps also fits into a factorization system,
and the projection $s:\Ar^{\pm \cap 1}(\eq) \to \eq$ is a fibration
by Example~\ref{facto.exa}. We let $\adj = \Ar^{\pm \cap
  1}(\eq)$. It is fibered over $\eq$, hence over $\Delta$, and it is
immediate to check that $\adj$ is a $2$-category, and $s$ then
becomes a $2$-functor $\nu$ of \eqref{nat.adj.eq}. The projection
$\nu$ has an obvious fully faithful left-adjoint $\eta:\eq \to \adj$
sending $\langle [n],e \rangle$ to the identity arrow $\id:\langle
[n],e \rangle \to \langle [n],e \rangle$. While it does not have a
right-adjoint, we can nevertheless define a $2$-subca\-tegory $\adj^0
\subset \adj$ by the cartesian square
\begin{equation}\label{adj.0}
\begin{CD}
\adj^0 @>{\beta}>> \adj\\
@V{\nu^0}VV @VV{\nu}V\\
\nat @>{\beta}>> \eq,
\end{CD}
\end{equation}
and $\nu^0$ has both a left-adjoint $\eta^0$ induced by $\eta$, and
a right-adjoint $\delta^0$. We can then define the other $2$-functor
in \eqref{nat.adj.eq} by $\delta = \beta \circ \delta^0$.

Explicitly, the $2$-category $\adj$ has two objects $0$, $1$, and by
definition, the endomorphism category $\adj(0,0)$ with its monoidal
structure is the category $\Delta_{\pm}$ with the monoidal structure
$B\Delta_{\pm} = \Ar^b(\Delta)$ of Subsection~\ref{mon.subs}. A
moment's reflection shows that $\adj(1,1)$ is the opposite category
$\Delta^<$ with the monoidal structure $B\Delta^< \cong
\Ar^c(\Delta) \subset \Ar^b(\Delta)$. Altogether, we have
\begin{equation}\label{2adj.01}
\adj(1,1) = \Delta^<, \ \adj(1,0)=\Delta_-, \
\adj(0,1)=\Delta_+, \ \adj(0,0) = \Delta_{\pm},
\end{equation}
with compositions
$$
\Delta^< \times \Delta_- \to \Delta_-, \ \Delta_+ \times
\Delta^< \to \Delta_+, \ \Delta^< \times \Delta^< \to \Delta^<,
\ \Delta_+ \times \Delta_- \to \Delta_{\pm}
$$
given by the concatenation product $- \circ -$, and compositions
$$
\Delta_{\pm} \times \Delta_+ \to \Delta_+, \ \Delta_- \times
\Delta_{\pm} \to \Delta_-, \ \Delta_{\pm} \times \Delta_{\pm} \to
\Delta_{\pm}, \ \Delta_- \times \Delta_+ \to \Delta^<
$$
given by the reduced concatenation product $- * -$. Then we have two
morphisms $\f:0 \to 1$, $\f^\vee:1 \to 0$ in $\adj$ corresponding to
the initial objects in $\Delta_+$, $\Delta_-$, and $\f^\vee \circ \f
\cong [0] \in \Delta^{<o} \cong \Delta_{\pm}$, $\f \circ \f^\vee
\cong [0] \in \Delta^<$, so that the unique map $o = \emptyset \to
      [0]$ produces maps $\aaa:\f^\vee \circ \f \to \id_0$,
      $\aaa^\vee:\id_1 \to \f \circ \f^\vee$.

\begin{lemma}\label{adj.le}
The quadruple $\langle \f,\f^\vee,\aaa,\aaa^\vee\rangle$ is an
adjoint pair in the $2$-category $\adj$, and for any adjoint pair
$\langle f,f^\vee,a,a^\vee \rangle$ in a $2$-category $\C$, there
exists a $2$-functor $\gamma:\adj \to \C$ and an isomorphism
$\gamma(\langle \f,\f^\vee,\aaa,\aaa^\vee\rangle) \cong \langle
f,f^\vee,a,a^\vee \rangle$, unique up to a unique isomorphism.
\end{lemma}

\proof{} Both compositions in \eqref{adj.rel} for the quadruple
$\langle \f,\f^\vee,\aaa,\aaa^\vee\rangle$ are given by the
composition $e \circ s:[0] \to [1] \to [0]$ in $\Delta_+ \cong
\Delta_-^o$, and since $e \circ s = \id$, we indeed have an adjoint
pair in $\adj$. Since all the categories in \eqref{2adj.01} have no
non-trivial isomorphisms, the $2$-category $\adj$ is strict, and the
underlying $1$-category $\rAdj$ is simply the path category of the
wheel quiver with two vertices $0$, $1$ and two edges $\f$,
$\f^\vee$. Therefore by \cite[Lemma 6.15]{ka.hh}, for any $\C$ and
adjoint pair $\langle f,f^\vee,a,a^\vee \rangle$, we have a unique
$2$-functor $\gamma$ from $\rAdj$ to $\C$ sending $\f$ to $f$ and
$\f^\vee$ to $f^\vee$, and we need to extend it to $\adj \supset
\Delta\rAdj$. In other words, we need to define $\gamma$ on all
the morphisms in the categories \eqref{2adj.01}. The morphisms in
$\Delta^< = \adj(1,1)$ are generated by the surjective degeneracy
maps $s^l_n:[n+1] \to [n]$ and injective face maps $d^l_n:[n-1] \to
[n]$, $l \in [n]$, and all the maps $\gamma(d^l_n)$, $\gamma(s^l_n)$
are uniquely defined by $\gamma(d^0_0)=\gamma(\aaa^\vee)=a^\vee$ and
$\gamma(s^0_0) = \gamma(\id_{\f} \circ \aaa \circ \id_{\f^\vee})=
\id_f \circ \operatorname{a} \circ \id_{f^\vee}$. The fact that
\eqref{adj.rel} then yields all the relations between the face and
degeneracy maps is well-known (when $\C=\cCat$, it just means that
$f \circ \f^\vee$ is a comonad, that is, a counital coassociative
coalgebra in $\C(1,1)$). This defines $\gamma$ on $\adj(1,1)$. For
the other three morphism categories \eqref{2adj.01}, note that they
admit faithful embeddings into $\adj(1,1)$ given by $\f \circ -$, $-
\circ \f^\vee$, $\f \circ - \circ \f^\vee$, and the morphisms that
are in the images of these embeddings are again generated by face
and degeneracy maps, modulo the same relations as in the ambient
category $\Delta^<$. Therefore $a$ and $a^\vee$ also uniquely define
$\gamma$ on $\adj(1,0),\adj(0,1),\adj(0,0) \subset \adj(1,1)$.
\endproof

\begin{remark}
Note that the Joyal duality \eqref{jo.eq} identifies the
$2$-category $\adj$ with its $2$-opposite $2$-category $\adj^\tau$
(the equivalence $\adj \cong \adj^\tau$ interchanges the objects $0$
and $1$). We of course also have $\eq \cong \eq^\tau$.
\end{remark}

\begin{remark}
One can also distinguish an intermediate class between adjoint pairs
and equivalences, by requiring that $a^\vee$ but not necessairly $a$
is invertible. In $\cCat$, this amounts to saying that the adjoint
$f^\vee$ to the reflexive functor $f$ is fully faithful. The
corresponding universal $2$-category is then the full subcategory
$\adj^p \subset \adj$ spanned by surjective bispecial $1$-special
maps. While the embedding $\adj^p \subset \adj$ is only a co-lax
$2$-functor, it has an honest right-adjoint $2$-functor $\adj \to
\adj^p$ that allows one to fit $\adj^p$ into the space in
\eqref{nat.adj.eq} between $\adj$ and $\eq$. The $2$-category
$\adj^p$ appeared in \cite[Subsection 8.5]{ka.hh} but only in
passing; it might be interesting to explore its combinatorics.
\end{remark}

\subsection{Co-adjoint pairs and twisting.}

Once we have constructed the universal $2$-categories
\eqref{nat.adj.eq}, we can define morphisms, adjoint pairs and
equivalences between $2$-functors in the same way as in
\eqref{ga.1}. Namely, assume given $2$-categories $\C$, $\C'$ and
$2$-functors $\gamma_0,\gamma_1:\C \to \C'$. Then a {\em
  $1$-morphism}, or a {\em natural transformation} from $\gamma_0$
to $\gamma_1$ is a $2$-functor
\begin{equation}\label{ga.nat}
\gamma:\C \times^2 \nat \to \C'
\end{equation}
that restricts to $\gamma_l$ on $\C \times \{l\} \subset \C \times^2
\nat$, $l = 0,1$, and similarly, an {\em adjoint pair} and a {\em
  pair of adjoint equivalences} are $2$-functors
\begin{equation}\label{ga.adj.eq}
\gamma:\C \times^2 \adj \to \C', \qquad \gamma:\C \times^2 \eq \to \C'.
\end{equation}
However, our explicit construction of the $2$-category $\adj$ also
reveals some additional structure that is not obvious from the
definitions. Namely, we also have a projection $t:\adj \cong
\Ar^{\pm \cap 1}(\eq) \to \eq \to \Delta$ sending an arrow to its
target, and for any $2$-category $\C$, we can define a category
$\C\{\adj\}$ by the cartesian square
\begin{equation}\label{C.adj.sq}
\begin{CD}
\C\{\adj\} @>>> \adj\\
@VVV @VV{t}V\\
\C @>>> \Delta.
\end{CD}
\end{equation}
The structural fibration $s:\adj \to \Delta$ then provides a
fibration $\C\{\adj\} \to \Delta$, and it is elementary to check
that this turns $\C\{\adj\}$ into a $2$-category.

\begin{defn}\label{coadj.def}
for any $2$-categories $\C$, $\C'$, a {\em coadjoint pair} of
$2$-functors from $\C$ to $\C'$ is a $2$-functor
$$
\gamma:\C\{\adj\} \to \C'.
$$
\end{defn}

If $\C \cong \C_{[0]} \times \Delta$ is discrete, then $\C\{\adj\}
\cong \C \times^2 \adj$, so that adjoint and co-adjoint pairs
defined on $\C$ coincide. It turns out that in general, co-adjoint
pairs are much easier to control. Namely, note that any $2$-functor
$\C \to \C'$ with discrete $\C'$ trivially satisfies the assumptions
on Lemma~\ref{fib.le}, so that it is a fibration. For any
$2$-category $\C$ equipped with a $2$-functor $\C \to \eq$, say that
a map in $\C$ is {\em $l$-special}, $l = 0,1$, if it is a cartesian
lifting of an $l$-special map in $\eq$. For any $2$-category $\C$,
denote $\C\{\eq\} = \C \times^2 \eq$, with its projection $\C\{\eq\}
\to \eq$, and say that a co-lax $2$-functor $\C \to \C'$ is {\em
  $l$-special} if it sends $l$-special maps in $\C\{\eq\}$ to
cartesian maps in $\C'$. Now, as we have mentioned, the functor $\nu$ of
\eqref{nat.adj.eq} has a left-adjoint $\eta:\eq \to \adj$, and for
any $\C$, these induce an adjoint pair of functors
\begin{equation}\label{nu.eta}
\nu:\C\{\adj\} \to \C\{\eq\}, \qquad \eta:\C\{\eq\} \to \C\{\adj\},
\end{equation}
where $\nu$ is a $2$-functor, and $\eta$ is a co-lax $2$-functor by
adjunction. We then have the following result.

\begin{lemma}[{{\cite[Lemma 7.4]{ka.hh}}}]\label{coadj.le}
For any $2$-category $\C$, the co-lax $2$-fun\-ctor $\eta$ of
\eqref{nu.eta} is $0$-special, and any $0$-special functor
$\gamma:\C\{\eq\} \to \C'$ uniquely factors as $\gamma \cong \gamma'
\circ \eta$ for a unique co-adjoint pair $\gamma':\C\{\adj\} \to
\C'$.\endproof
\end{lemma}

Lemma~\ref{coadj.le} essentially says that a co-adjoint pair is
completely defined by the corresponding $0$-special co-lax
$2$-functor, and there is also one general result proved in
\cite[Subsection 7.4]{ka.hh} that helps to construct such
$2$-functors. To state it, we need some preliminaries. For any
co-lax $2$-functor $\gamma:\C' \to \C$, one can define a
$2$-category $\gamma^*\C$ by the cartesian square
$$
\begin{CD}
\gamma^*\C @>>> \C\\
@VVV @VVV\\
\eps_*\C'_{[0]} @>>> \eps_*\C_{[0]},
\end{CD}
$$
where as in \eqref{red.sq}, $\eps:\ppt \to \Delta$ is the embedding
onto $[0] \in \Delta$. Informally, $\gamma^*$ has the same objects
as $\C'$ and the same morphisms as $\C$. Then $\gamma$ factors as
\begin{equation}\label{ga.fac}
\begin{CD}
\C' @>{\wgamma}>> \gamma^*\C @>{\ogamma}>> \C,
\end{CD}
\end{equation}
and in keeping with our usage for the usual categories, we say that
$\gamma$ is {\em dense} resp.\ {\em full} if $\ogamma$
resp.\ $\wgamma$ is an equivalence. A $2$-subcategory $\C' \subset
\C$ is full resp.\ dense iff so is the embedding $2$-functor $\C'
\to \C$. For example, if $\gamma:\C'_{[0]} \times \Delta \to
\C_{[0]} \times \Delta \to \C$ is the embedding $2$-functor for some
$\C'_{[0]} \subset \C_{[0]}$, then $\gamma^*\C \subset \C$ is the
full $2$-subcategory spanned by objects $c \in \C'_{[0]}$.

\begin{prop}[{{\cite[Proposition 7.9]{ka.hh}}}]\label{tw.prop}
Assume given some $2$-cate\-go\-ries $\C'$, $\C$, denote by
$\iota_1:\C' \to \C'\{\eq\}$ the embedding onto $\C' \times \{1\}
\subset \C'\{\eq\}$, and assume given a $0$-special co-lax
$2$-functor $\gamma:\C'\{\eq\} \to \C$. Then the component
$\wgamma$ of its decomposition \eqref{ga.fac} factors as
\begin{equation}\label{dm.facto}
\begin{CD}
\C'\{\eq\} @>{\wgamma_1 \times \id}>> \gamma_1^*\C\{\eq\}
@>{\gamma^\dm}>> \gamma^*\C,
\end{CD}
\end{equation}
where $\gamma_1 = \gamma \circ \iota_1$, $\wgamma_1:\C' \to
\gamma_1^*\C$ is as in \eqref{ga.fac}, and $\gamma^\dm$ is a
$0$-special co-lax $2$-functor over $\eq/\Delta^o$ equipped with an
isomorphism $\iota_1^*(\gamma^\dm) \cong \id$. Moreover, such a
factorization is unique up to a unique isomorphism.\endproof
\end{prop}

In particular, Proposition~\ref{tw.prop} provides the ``twisting
co-lax $2$-functor'' $\Theta$ from $\gamma_1^*\C$ to $\C$ obtained
by restricting $\gamma^\dm$ to $\gamma_1^*\C \times \{0\} \subset
\gamma_1^*\C\{\eq\}$. Explicitly, for any $c \in \C'_{[0]}$,
$\gamma$ gives rise to an adjoint pair $\gamma_c:\adj \to \C$
consisting of objects $\gamma_0(c),\gamma_1(c) \in \C_{[0]}$,
morphisms $h_c \in \C(\gamma_0(c),\gamma_1(c))$, $h^\vee_c \in
\C(\gamma_1(c),\gamma_0(c))$, and adjunction maps between their
compositions. Then $\Theta$ sends $c$ to $\gamma_0(c)$, and on
morphisms, it is given by
\begin{equation}\label{theta.eq}
\Theta(g) = h^\vee_{c'} \circ g \circ h_c, \qquad g \in
(\gamma_1^*\C)(c,c') \cong \C(\gamma_1(c),\gamma_1(c')).
\end{equation}
The maps \eqref{fib.fu} are induced by the adjunction maps between
$h_\idot$ and $h_\idot^\vee$.

\begin{exa}\label{equi.exa}
For a somewhat trivial but useful application of
Proposition~\ref{tw.prop}, assume given a $2$-category $\C$ with a
full $2$-subcategory $\C' \subset \C$, and assume that for any $c
\in \C'_{[0]} \subset \C_{[0]}$, we are given an equivalent object
$\theta(c) \in \C_{[0]}$ --- that is, a $2$-functor $\gamma_c:\eq
\to \C$ sending $1$ to $c$ and $0$ to some $\theta(c)$. Then the
functors $\gamma_c$ together define a $2$-functor $\gamma:\C'_{[0]}
\times \eq \to \C$ such that $\gamma_1^*\C \cong \C'$, and then
Proposition~\ref{tw.prop} provides a co-lax $2$-functor $\Theta:\C'
\to \C$, $c \mapsto \theta(c)$. By \eqref{theta.eq}, since we are
dealing with equivalences and not just adjoint pairs, $\Theta$ is in
fact a $2$-functor and a full embedding $\C' \to \C$ different from
the original one. However, the two embeddings are equivalent in the
sense of \eqref{ga.adj.eq}, with the equivalence provided by
$\gamma^\dm$ (that in this case is also a $2$-functor). Note that
the equivalences $\gamma_c$ are chosen separately and independently
for each $c$, and do not need to be compatible with morphisms in
$\C'$ in any way.
\end{exa}

\subsection{Iterated adjoint pairs.}

Let us now describe a non-trivial application of
Proposition~\ref{tw.prop} given in \cite[Subsection 7.5]{ka.hh}. It
concerns compositions of reflexive morphisms and adjoint pairs.

For any objects $c,c' \in \C_{[0]}$ in a $2$-category $\C$, adjoint
pairs of maps between $c$ to $c'$ form a category $\aAdj(\C)(c,c')$.
In fact, the category is a groupoid, and the forgetful functor
$\aAdj(c,c') \to \C(c,c')_{\Iso}$, $\langle f, f^\vee,a,a^\vee
\rangle \mapsto f$ is fully faithful, so we may identify
$\aAdj(\C)(c,c')$ with its essential image and treat is a
subcategory in $\C(c,c')$. This is a stronger form of the usual
uniqueness of adjoints; for $\C = \cCat$, it is easy to check it
directly, and the general case again reduces to $\cCat$ by the
Yoneda embedding of Theorem~\ref{rigi.thm}. Equivalences then span a
full subcategory $\eEq(\C)(c,c') \subset \aAdj(\C)(c,c')$. Since
both reflexive morphisms and equivalences are closed under
compositions, we in fact have dense $2$-subcategories
\begin{equation}\label{eeq.aadj}
\eEq(\C) \subset \aAdj(\C) \subset \C
\end{equation}
with morphism categories $\eEq(\C)(-,-)$, $\aAdj(\C)(-,-)$, and
moreover, by Lemma~\ref{adj.le}, we have $\aAdj(\C)_{[1]} \cong
\Fun^2(\adj,\C)$ and $\eEq(\C)_{[1]} \cong \Fun^2(\eq,\C)$. Our goal
is a similar universal description for the whole $2$-categories
\eqref{eeq.aadj}.

To package the answer, it is convenient to say that a {\em
  $2$-kernel} is a small category $\K$ equipped with a fibration $\K
\to \Delta^o \times \Delta$ such that for any $[n] \in \Delta^o$,
$\K_{[n]} \to \Delta$ is a $2$-category (this is a version of
``$\Delta$-kernels'' used in \cite[Subsection 7.5]{ka.hh}). For any
$2$-kernel $\K$ and $2$-category $\C$, we then have a fibration
$\Fun^2(\K,\C) \to \Delta$ with fibers $\Fun^2(\K,\C)_{[n]} \cong
\Fun^2(\K_{[n]},\C)$ and transition functors given by pullbacks with
respect to the transition functors $\K_{[n]} \to \K_{[m]}$ of the
fibration $\K \to \Delta^o$. We then want to promote
\eqref{nat.adj.eq} to a sequence
\begin{equation}\label{Nat.Adj.Eq}
\begin{CD}
\Nat @>{\delta}>> \Adj @>{\nu}>> \Eq
\end{CD}
\end{equation}
of $2$-kernel and functors cartesian over $\Delta^o \times \Delta$
that produces \eqref{eeq.aadj} after applying $\Fun^2(-,\C)$.

For $\Nat$, the answer is tautological: we have $\Nat_{[n]}
\cong\Delta [n]$ for any $[n]$, and $\Nat$ itself is the twisted
arrow category $\Tw(\Delta^o)$ of Example~\ref{tw.exa}. For $\Eq$,
it is easy to see that we have $\Eq_{[n]} = \Delta
e(\{0,\dots,n\})$, and if we let $V:\Delta \to \Gamma$ be the
functor sending an ordinal $[n]$ to the finite set
$V([n])=\{0,\dots,n\}$, then $\Eq \cong (V^o \times
V)^*\Tw(\Gamma^o)$, and $V$ induces a cartesian functor $\beta:\Nat
\to \Eq$.

To construct $\Adj$, consider the cofibration $t:\Ar^t(\Delta) \cong
\Delta^\hdot \to \Delta$ of Remark~\ref{del.dot.rem}, with the
functor $s:\Delta^\hdot \to \Delta$, $\langle [n],l \rangle \mapsto
[n-l]$, and note that if we restrict $s$ to $\theta^*\Delta^\hdot
\subset \Delta^\hdot$, where the embedding $\theta:\Delta^o \subset
\Delta^{<o} \cong \Delta_{\pm} \subset \Delta$ is induced by the
Joyal duality \eqref{jo.eq}, then it factors through $\Delta_-
\subset \Delta$. We can then apply the duality \eqref{p.m.o} and
obtain a functor $w':\theta^*\Delta^\hdot \to \Delta_+^o \subset
\Delta^o$. Moreover, we can modify it by setting
\begin{equation}\label{w.eq}
w(\langle [n],l \rangle) = \begin{cases} [n-l] = w'(\langle [n],l
    \rangle), &\quad l > 0,\\
t([n-1]) \subset [n] = w'(\langle [n],l), &\quad l =0,
\end{cases}
\end{equation}
and we note that this gives a well-defined functor
$w:\theta^*\Delta^\hdot \to \Delta^o$. Indeed, for any map
$f:\langle [n],l \rangle \to \langle n',l' \rangle$ in
$\theta^*\Delta^\hdot$, the map $f$ must be bispecial, so that if
$l=0$ then $l'=0$, and moreover, $w'(f) = f^\dg:[n'] \to [n]$ sends
$t([n-1]) \subset [n]$ into $t([n' - 1] \subset [n']$.  Therefore
for any $f$, $w'(f):[n'-l'] \to [n-l]$ sends $w(\langle [n'],l'
\rangle) \subset w'(\langle [n'],l' \rangle)$ into $w(\langle [n],l
\rangle) \subset w'(\langle [n],l \rangle)$, and the functor $w$
given by \eqref{w.eq} is indeed well-defined. Explicitly, for any
$[n] \in \Delta^o$, it restricts to a functor $w_n:[n+1]^o =
(\theta^*\Delta^\hdot)_{[n]} \to \Delta^o$, and then the opposite
functor $w_n^o:[n+1] \to \Delta$ is given by the diagram
\begin{equation}\label{n.dia}
\begin{CD}
[n] @>{g_0}>> [n] @>{g_1}>> [n-1] @>{g_2}>> \dots @>{g_{n-1}}>> [1]
@>{g_n}>> [0]
\end{CD}
\end{equation}
in $\Delta$, where $g_0 = \id$, and $g_l:[n+1-l] \to [n-l]$ for $l
\geq 1$ sends $0$ to $0$ and $i \in [n+1-l]$, $i \geq 1$ to $i-1 \in
     [n-l]$.

Now consider the fibration $\Eq \to \Delta^o$, and let $\Adj'' =
(\theta^*t)_{**}w^*\Eq$, where $\theta^*t:\theta^*\Delta^\hdot \to
\Delta^o$ is the restriction of the cofibration $t:\Delta^\hdot \to
\Delta$, and $(\theta^*t)_{**}$ is as in \eqref{st.st}. By
definition, $\Adj''$ comes equipped with a fibration $\pi:\Adj'' \to
\Delta^o$, and the fiber $\Adj''_{[n]}$ over some $[n] \in \Delta^o$
is the category of pairs $\langle m_\idot,\phi \rangle$ of a functor
$m_\idot:\theta([n]) = [n+1] \to \Delta$ and a map $\phi:V \circ
m_\idot \to V \circ w^o_n$. Explicitly, $m_\idot$ is a diagram
\begin{equation}\label{m.dia}
\begin{CD}
[m_0] @>{f_0}>> [m_1] @>{f_1}>> \dots @>{f_n}>> [m_{n+1}],
\end{CD}
\end{equation}
and for any $l \in [n+1]$, $\phi$ provides a partition $[m_l] =
[m_l]_0 \copr \dots \copr [m_l]_j$ into $j$ subsets, where $j =
\min(n+1-l,n)$. If $l>0$, then the map $f_l$ in \eqref{m.dia} sends
$[m_l]_{i+1}$ into $[m_{l+1}]_i$, $0 \leq i < j$, and it sends
$[m_l]_0$ into $[m_{l+1}]_0$. The map $f_0$ simply sends $[m_0]_i$
into $[m_1]_i$.

We now observe that since $\theta:\Delta^o \to \Delta$ sends all
maps to bispecial maps, both $[m_0]$ and $m_{[n+1]}$ in
\eqref{m.dia} are functorial with respect to $\langle m_\idot,\phi
\rangle$, so that we have two projections $\sigma,\tau:\Adj'' \to
\Delta$. We let $\Adj' = \tau^{-1}([0]) \subset \Adj''$ be the full
subcategory spanned by pairs $\langle m_\idot,\phi \rangle$ with
$[m_{n+1}] = [0]$, and we still have the projection $\sigma:\Adj'
\to \Delta$ and the fibration $\pi:\Adj' \to \Delta^o$. Moreover,
say that a pair $\langle m_\idot,\phi \rangle$ is {\em admissible}
if $f_l:[m_l]_i \to [m_{l+1}]_{g_l(i)}$ is an isomorphism whenever
$g_l(i) \geq 1$ (note that by \eqref{n.dia}, we have $g_l(i) = i-1$
if $l \geq 1$ and $g_l(i) = i$ if $l=0$). Let $\Adj \subset \Adj'$
be the full subcategory of admissible pairs.

\begin{lemma}
The product $\pi \times \sigma:\Adj' \to \Delta^o \times \Delta$ is
a $2$-kernel, and so is the induced projection $\Adj \subset \Adj'
\to \Delta^o \times \Delta$. Moreover, the embedding $\Adj \subset
\Adj'$ admits a right-adjoint functor $\Adj' \to \Adj$ cartesian
over $\Delta^o \times \Delta$.
\end{lemma}

\proof{} This is essentially \cite[Lemma 7.10]{ka.hh} (modulo the
change of variance explained in Remark~\ref{cova.rem}).
\endproof

One further observes that $[m_0]$ in a diagram \eqref{m.dia} comes
equipped with a partition into $n$ subsets, thus defines an object
in $\Eq_{[n]}$, and the fibration $\Adj \to \Delta^o \times \Delta$
then factors through a functor $\nu:\Adj \to \Eq$ cartesian over
$\Delta^o \times \Delta$. We can then define the subcategory $\Adj^0
\subset \Adj$ by the cartesian square
\begin{equation}\label{Adj.0}
\begin{CD}
\Adj^0 @>{\beta}>> \Adj\\
@V{\nu^0}VV @VV{\nu}V\\
\Nat @>{\beta}>> \Eq,
\end{CD}
\end{equation}
the universal version of \eqref{adj.0}, and $\nu^0$ again has a
right-adjoint $\delta^0$, so we can let $\delta = \beta \circ
\delta^0$ and obtain our diagram \eqref{Nat.Adj.Eq}.

\begin{exa}\label{ite.exa}
If $n=0$, then an admissible pair $\langle m_\idot,\phi \rangle$
consists of a diagram $[m] \to [0]$; since $w_0$ is the diagram $[0]
\to [0]$, $\phi$ is unique. Therefore we have $\Adj_{[0]} = \Delta$.

If $n=1$, then we have a diagram $[m_0] \to [m_1] \to [0]$, and
$\phi$ defines a partition $[m_1] = [m_1]_0 \copr [m_1]_1$ and a
partition $[m_0]=[m_0]_0 \copr [m_0]_1$. The latter is uniquely
defined by the former, and what we have is actually a map $f_0$ in
$\eq$. Admissibility means that the map $1$-special, so that
$\Adj_{[1]} \cong \adj$.

We have two maps $s,t:[0] \to [1]$ in $\Delta$, and the bispecial
maps $p = \theta(s^o)$ resp.\ $q = \theta(t^o)$ from $[2]$ to $[1]$
send $1 \in [2]$ to $0$ resp.\ $1$. The transition functor
$s^{o*}:\Delta \to \adj$ of the fibration $\Adj \to \Delta^o$ then
sends a diagram $[m] \to [0]$ to its pullback $p^*$ with respect to
$p$ --- that is, to the diagram $[m] \to [m] \to [0]$, with the
first map being the identity map --- while $t^{o*}$ sends it to $[m]
\to [0] \to [0]$, the pullback $q^*$ with respect to $q$. The
partition data are prescribed by the adjoint maps $p^\dg,q^\dg:[1]
\to [2]$, and we have $p^\dg(1)=2$ and $q^\dg(1)=1$. Therefore
$w'(p):[1] \to [2]$ sends $1$ to $2$, so that $w(p):[0] \to [1]$
sends $0$ to $1$, while $w(q)$ sends $0$ to $0$. This means that for
the diagram $[m] \to [m] \to [0]$, we take $[m] = [m]_1$, while for
the diagram $[m] \to [0] \to [0]$, we take $[m] = [m]_0$. Both
diagrams are therefore admissible, and $s^{o*}$, $t^{o*}$ are the
embeddings $\Delta \to \adj$ onto $1$ and $0$.
\end{exa}

\begin{prop}[{{\cite[Proposition 7.12]{ka.hh}}}]\label{adj.prop}
For any $2$-category $\C$, the fibration $\Fun^2(\Adj,\C) \to
\Delta$ is a $2$-category, and the functor
$$
\delta^*:\Fun^2(\Adj,\C) \to \Fun^2(\Nat,\C) \cong \C
$$
induced by \eqref{Nat.Adj.Eq} factors through an equivalence
$\Fun^2(\Adj,\C) \cong \aAdj(\C)$.\endproof
\end{prop}

Let us explain the main idea behind the proof of
Proposition~\ref{adj.prop}. It uses an alternative inductive
construction of the $2$-categories $\Adj_n = \Adj_{[n]}$ that also
explains how can one possibly come up with such a weird object.

Take some $n \geq 2$, consider the cocartesian square \eqref{seg.sq}
with $l=1$, and let $e:[n] \to [1] \times [n-1]$ be the map given by
the embedding onto $[1] \times \{0\}$ on $s([1]) \subset [n]$, and
by the embedding onto $\{1\} \times [n-1]$ on $t([n-1]) \subset
[n]$. Note that $e$ is a map of partially ordered sets, thus a
functor, and it has a left-adjoint functor $p:[1] \times [n-1] \to
[n]$ sending $i \times j$ to $i$ if $j=0$ and to $j+1$ if $j \geq
1$. Consider the $2$-category $\Adj_{n-1}$ and the corresponding
$2$-category $\Adj_{n-1}\{\adj\}$ of \eqref{C.adj.sq}. Then by
definition, an object in $\Adj_{n-1}\{\adj\}$ is a diagram
\eqref{m.dia} equpped with two types of partition data. Firstly,
$f_0:[m_0] \to [m_1]$ defines an object in $\adj$, so that we have a
map $\phi_1:V([m_1]) \to \{0,1\}$ (and the induced map $\phi_1 \circ
V(f_0):V([m_0]) \to \{0,1\}$). Secondly, the part of the diagram
that starts from $[m_1]$ defines an object in $\Adj_{n-1}$, so that
we have a map $\phi_{n-1}:V([m_1]) \to \{0,\dots,n-1\}$. Altogether,
we have a map $\phi_1 \times \phi_{n-1}:V([m_1]) \to V([1] \times
[n-1]) = V([1]) \times V([n-1])$ and the induced map $(\phi_1 \times
\phi_{n-1}) \circ V(f_0):V([m_0]) \to V([1] \times [n-1])$.

Observe that if the image of the map $\phi_1 \times \phi_{n-1}$ lies
inside the image of the map $V(e):V([n]) \to V([1] \times [n-1])$,
then we can extend $\phi_{n-1}$ to a map $\phi_n:V(m_\idot) \to
V(w_n^o)$ by setting $\phi_n = (\phi_1 \times \phi_{n-1}) \circ
V(f_0)$ on $V([m_0])$, $\phi_n = \phi_1 \times \phi_{n-1}$ on
$V([m_1])$, and $\phi_n = \phi_{n-1}$ on $V([m_l])$ for $l \geq
2$. One checks immediately that the admissibility conditions match,
so that we obtain a cartesian square
$$
\begin{CD}
\Adj_n @>{e}>> \Adj_{n-1}\{\adj\}\\
@VVV @VVV\\
\Eq_n @>{e}>> \Eq_{n-1}\{\eq\},
\end{CD}
$$
where we identify $\Eq_{n-1}\{\eq\} \cong \Delta e(V([1] \times
[n]))$. Moreover, the map $p$ extends to a $2$-functor
$p:\Adj_{n-1}\{\adj\} \to \Adj_n$ such that $p \circ e \cong
\Id$. Indeed, by Lemma~\ref{coadj.le} and Proposition~\ref{tw.prop},
such an extension is uniquely defined by its compositions $p_i = p
\circ a_i$ with the embeddings $a_i:\adj \subset \Adj_{n-1}\{\adj\}$
onto $i \times \adj$, $0 \leq i \leq n-1$, and we can take $p_0$ to
be the transition functor $s^{o*}:\adj = \Adj_1\to \Adj_n$ of the
fibration $\Adj \to \Delta^o$ corresponding to the embedding $s:[1]
\to [n]$, and let $p_i$, $i \geq 1$ correspond to the projection $[1]
\to [0] \to [n]$ onto $i+1 \in [n]$. We conclude that $\Adj_n$ is a
retract of the $2$-category $\Adj_{n-1}\{\adj\}$.

Now, by Example~\ref{ite.exa}, we know that $\aAdj(\C) \cong
\Fun^2(\Adj,\C)$ over the objects $[0]$ and $[1]$ in $\Delta$, so to
prove Proposition~\ref{adj.prop}, it suffices to check that
$\Fun^2(\Adj,\C)$ satisfies the Segal condition. Moreover, by
induction, it suffices to check it for a square \eqref{seg.sq} with
$l=1$, assuming that we know that $\aAdj(\C)_{[i]} \cong
\Fun^2(\Adj_i,\C)$ for $i \leq n-1$. What we have to do, then, is to
take a $2$-functor $\gamma_1:\adj \to \C$ and a $2$-functor
$\gamma_{n-1}:\Adj_{n-1} \to \C$ such that
$\gamma_1(1)=\gamma_{n-1}(0)$ is the same object $c \in \C_{[0]}$,
and show that there exists a unique $2$-functor $\gamma_n:\Adj_n \to
\C$ equipped with isomorphisms $\gamma_n \circ s^{o*} \cong
\gamma_1$ and $\gamma_n \circ t^{o*} \cong \gamma_{n-1}$. However,
again by Lemma~\ref{coadj.le} and Proposition~\ref{tw.prop}, for any
$(n-1)$-tuple of adjoint pairs $\gamma'_i:\adj \to \C$ sending $1
\in \adj_{[1]}$ to the object $\gamma_{n-1}(i) \in \C_{[0]}$, the
$2$-functor $\gamma_{n-1}$ extends uniquely to a co-adjoint pair
$\gamma:\Adj_{n-1}\{\adj\} \to \C$ equipped with isomorphisms
$\gamma \circ a_i \cong \gamma'_i$. It remains to observe that
$\gamma$ factors through the retract $\Adj_n$ of the $2$-category
$\Adj_{n-1}\{\adj\}$ if and only if for any $i \geq 1$, $\gamma'_i$
factors through the embedding $\Delta \to \C$ onto
$\gamma_{n-1}(i)$, and take $\gamma'_0 = \gamma_1$.

\section{Applications.}\label{app.sec}

\subsection{Symmetric monoidal structures.}\label{comm.subs}

In order to show how the formalism we have developed applies to real
life, it is convenient to start with yet another piece of formalism
--- namely, with a description of symmetric monoidal
categories. Traditionally, these are defined in terms of
associativity and commutativity isomorphisms that satisfy higher
constraints (the pentagon and the hexagon axiom, and the unitality
axioms that are often ignored). This theory has no strict version,
since the commutativity isomorphism is almost never an identity.

For a description of symmetric monoidal categories in the spirit of
Definition~\ref{2cat.def}, let $\Gamma$ be the category of finite
sets, and let $\Gamma_+$ be the category of finite sets and
partially defined maps -- that is, maps from $S_0$ to $S_1$ are
isomorphism classes of diagrams
\begin{equation}\label{dom}
\begin{CD}
S @<{i}<< \wt{S} @>{f}>> S'
\end{CD}
\end{equation}
in $\Gamma$ with injective $i$, with compositions given by fibered
products. Equivalently, $\Gamma_+$ is the category of finite pointed
sets, with the equivalence sending a set $S_+$ with the
distinguished element $o \in S_+$ to the complement $S = S_+ \setminus
\{o\} \subset S_+$, and a map $f:S_+ \to S'_+$ to the diagram
\eqref{dom} with $\wt{S} = f^{-1}(S') \subset S$. A map \eqref{dom}
in $\Gamma_+$ is {\em anchor} resp.\ {\em structural} if $f$
resp.\ $i$ is invertible. Coproducts in $\Gamma$ are also coproducts in
$\Gamma_+$, and for any $S_0,S_1 \in \Gamma_+$, we have anchor maps
\begin{equation}\label{anch.eq}
a_0:S_0 \copr S_1 \to S_0, \qquad a_1:S_1 \copr S_1 \to S_1
\end{equation}
defined by the embeddings $i_0:S_0 \to S_0 \copr S_1$, $i_1:S_1 \to
S_0 \copr S_1$.

\begin{defn}\label{sym.def}
A fibration $\C \to \Gamma_+^o$ {\em satisfies the Segal condition}
if $\C_\emptyset$ is equivalent to the point category $\ppt$, and
for any $S_0,S_1 \in \Gamma$, the projection
$$
a_0^* \times a_1^*:\C_{S_0 \copr S_1} \to \C_{S_0} \times \C_{S_1}
$$
induced by the anchor maps \eqref{anch.eq} is an equivalence of
categories. A {\em unital symmetric monoidal structure} on a
category $\C$ is a fibration $\Bi\C \to \Gamma^o_+$ that satisfies
the Segal condition and is equipped with an equivalence
$\Bi\C_{\ppt} \cong \C$. The {\em opposite unital symmetric monoidal
  structure} on the opposite category $\C^o$ is given by $\Bi\C^o =
(\Bi\C)_\perp^o$.  A monoidal structure on a functor $\gamma:\C \to
\C'$ between two categories $\C$, $\C'$ equipped with unital
symmetric monoidal structures $\Bi\C$, $\Bi\C'$ is given by a
functor $\Bi\gamma:\Bi\C \to \Bi\C'$, cartesian over $\Gamma^o_+$
and equipped with an isomorphism $\Bi\gamma(\ppt) \cong \gamma$. A
{\em co-lax symmetric monoidal structure} on $\gamma$ is given by a
functor $\Bi\gamma:\Bi\C \to \Bi\C'$ that is cartesian over the
anchor maps, and a lax symmetric monoidal structure on $\gamma$ is a
co-lax one on the opposite functor $\gamma^o$.
\end{defn}

We note that this is quite parallel to Definition~\ref{mon.def} and
Definition~\ref{lax.def}. In fact, the standard {\em simplicial
  circle} $\Sigma:\Delta^o \to \Sets$ is obtained by gluing together
the two ends of the standard $1$-simplex, $\Sigma([n])$ is finite
for any $[n]$ and has the distinguished point given by glued ends,
and the opposite functor $\Sigma^o:\Delta \to \Gamma^o_+$ sends the
maps $s$ and $t$ in \eqref{seg.eq} to the maps \eqref{anch.eq} in
$\Gamma_+^o$ (and more generally, a map $f$ is an anchor map in
$\Delta$ if and only if $\Sigma^o(f)$ is an anchor map in
$\Gamma^o_+$). Therefore for any unital symmetric monoidal structure
$\Bi\C$ on a category $\C$, the pullback $B\C = \Sigma^{o*}\Bi\C$ is
a unital monoidal structure. A symmetric monoidal structure on a
functor restricts to a monoidal structure in the sense of
Definition~\ref{mon.def}, and a lax or a co-lax one restricts to a
corresponding structure in the sense of
Definition~\ref{lax.def}. Just as in the non-symmetric case, if a
functor $\gamma:\C \to \C'$ admits a right-adjoint $\gamma^\dg$,
then a co-lax symmetric monoidal structure on $\gamma$ induces a lax
one on $\gamma^\dg$, and vice versa.

Explicitly, for any unital symmetric monoidal structure $\Bi\C$, the
anchor maps in $\Gamma_+$ serve to provide identification $\Bi\C_S
\cong \C^S$, $S \in \Gamma$, thus anchoring the {\em a priori}
arbitrary fibers $\Bi\C_S$ of the fibration $\Bi\C$ to our category
$\C$. Structural maps then encode the structure: the tensor product
functor $\C \times \C \to \C$ is induced by the unique structural map $\ppt
\copr \ppt \to \ppt$, the unit object corresponds to the map
$\emptyset \to \ppt$, and various associativity and commutativity
isomorphisms are packaged into the isomorphisms \eqref{fib.eq} for
the fibration $\Bi\C \to \Gamma_+$. This can be made into an
explicit comparison theorem but we will not need it.

\begin{remark}
The Segal condition of Definition~\ref{sym.def} is the original
condition introduced by Segal in \cite{seg} (although he worked with
topological spaces and weak equivalences rather then categories and
equivalences). The (obvious) generalization to $2$-categories
appeared slightly later and inherited the name.
\end{remark}

The first example of a unital symmetric monoidal structure is given
by cartesian products: if a category $\C$ has finite products, they
define a unital symmetric monoidal structure $\Bi\C$. Explicitly,
objects of $\Bi\C$ are pairs $\langle S,c_\idot \rangle$ of a finite
set $S \in \Gamma$ and a collection of objects $c_s \in \C$, $s \in
S$. A morphism $\langle S',c'_\idot \rangle \to \langle S,c_\idot
\rangle$ is a pair of a diagram \eqref{dom} and a collection of maps
$c'_{f(s)} \to c_s$, $s \in \wt{S} \subset S$. Note that the
embedding $i_c:\ppt \to \C$ onto any object $c \in \C$ carries a
unique co-lax unital symmetric monodial structure $\Bi i_c$ that
sends a set $S$ to the collection $c_s=c$ of copies of the object
$c$, and a morphism $f$ to the collection of identity maps $\id:c
\to c$. One can also consider the transpose cofibration
$(\Bi\C)_\perp \to \Gamma_+$. It has the same objects, and a
morphism $\langle S,c_\idot \rangle \to \langle S',c'_\idot \rangle$
is a diagram \eqref{dom} and a collection of maps
\begin{equation}\label{cop.eq}
\prod_{s \in f^{-1}(s')}c_s \to c_{s'}, \qquad s' \in S',
\end{equation}
where the product is the cartesian product in $\C$. Dually, if $\C$
has finite coproducts, they also define a unital symmetric monoidal
structure. In fact, a unital symmetric monoidal structure $\Bi\C$
on a category $\C$ defines a unital symmetric monoidal structure
$\Bi\C^o = (\Bi\C)_\perp^o$ on the opposite category $\C^o$, and if
$\C$ has finite coproducts, $\C^o$ has finite products.

The description of the cartesian product monoidal structure in terms
of the maps \eqref{cop.eq} may look slightly artificial but it has
its uses. For example, for any commutative ring $k$, the category
$k\amod$ is equipped with a forgetful functor $k\amod \to
\Sets$. One can then define a unital symmetric monoidal structure
$\Bi k\amod$ as follows: objects are pairs $\langle S,c_\idot
\rangle$, $S \in \Gamma$, $c_\idot \in k\amod$, and a morphism in
the transpose cofibration $(\Bi k\amod)_\perp$ is a diagram
\eqref{dom} and a collection of set-theoretic maps \eqref{cop.eq}
that are $k$-linear in each argument. This condition is obviously
closed under compositions, so we indeed obtain a well-defined
category. The corresponding product on $k\amod$ is the usual tensor
product $- \otimes_k -$ (in fact, this is how the tensor product of
vector spaces is defined in any good linear algebra textbook).

To extend this further, assume that $k$ is Noetherian. Then the
category $\Comm(k)$ of finitely generated commutative $k$-algebra
has finite coproducts, hence the coproduct monoidal structure
$\Bi\Comm(k)$, and the transpose cofibration $(\Bi\Comm(k))_\perp
\to \Gamma_+$ can be also described by imposing conditions on the
maps \eqref{cop.eq}. Moreover, the same procedure then constructs a
unital symmetric monoidal structure on the category $\Comm(k)\amod$
of Example~\ref{sch.exa} and on the forgetful functor
$\phi:\Comm(k)\amod \to \Comm(k)$. The product is given by $\langle
A,M \rangle \otimes \langle A',M' \rangle = \langle A \otimes_k A',M
\otimes_k M' \rangle$. The opposite categories $\Comm(k)^o$,
$\Comm(k)\amod^o$ and the opposite functor $\phi^o$ then also carry
unital symmetric monoidal structures, and moreover, the functor
$\Bi\phi^o:\Bi\Comm(k)\amod^o \to \Bi\Comm(k)^o$ is a
fibration. Furthermore, the category $\Sch(k)$ has the unital
symmetric monoidal structure given by cartesian products, the full
embedding $j:\Comm(k)^o \to \Sch(k)$ carries an obvious unital
symmetric monoidal structure $\Bi j$, and the pushforward fibration
$\Bi\QCoh(\Sch(k))^o = (\Bi j)_*\Bi\Comm(k)\amod^o$ defines a unital
symmetric monoidal structure on $\QCoh(\Sch(k))^o$, hence also on
$\QCoh(\Sch(k))$, with the product
\begin{equation}\label{box.eq}
\langle X,\F \rangle \otimes \langle X',\F' \rangle = \langle X
\times_{\Spec k} X',\F \boxtimes_k \F'\rangle.
\end{equation}
As in Example~\ref{sch.exa}, we can use the same procedure to obtain
a unital symmetric monoidal structure $\Bi\QCoh_\idot(\Sch(k))^o$ on
$\QCoh_\idot(\Sch(k))^o$, and since the tensor product of $h$-flat
complexes is $h$-flat, it induces a unital symmetric monoidal
structure $\Bi\QCoh_\idot^\flat(\Sch(k))^o$ on the full subcategory
$\QCoh_\idot^\flat(\Sch(k))^o \subset
\QCoh_\idot(\Sch(k))^o$. Localizing with respect to
quasiisomorphisms then gives a unital symmetric monoidal structure
$\Bi\D(\Sch(k))^o$ on $\D(\Sch(k))^o$ and the opposite symmetric
monoidal structure $\Bi\D(\Sch(k))$ on $\D(\Sch(k))$, with the
product given by the derived version of \eqref{box.eq}. The
forgetful bifibrations $\D(\Sch(k)) \to \Sch(k)^o$, $\D(\Sch(k))^o
\to \Sch(k)$ carry unital symmetric monoidal structures. Moreover,
for any $X \in \Sch(k)$, we can equip the embedding $i_X:\ppt \to
\Sch(k)$ with its unique co-lax unital symmetric monoidal structure
$\Bi i_X$, and then $\Bi \D(X)^o = (\Bi i_X)^*\Bi\D(\Sch(k))$
provides the usual symmetric monoidal structure on $\D(X)^o \cong
i_X^*\D(\Sch(k))^o$.

\subsection{Morita $2$-category.}

The first real-life example of a $2$-category that is not manifestly
strict is probably the Morita $2$-category: its objects are algebras
$A$ over a fixed commutative ring $k$, morphisms from $A_0$ to $A_1$
are $A_0\text{-}A_1$-bimodules (that is, $A_0^o \otimes_k
A_1$-modules), and the composition is given by the tensor
product. To construct it in terms of Definition~\ref{2cat.def}, one
can follow the procedure described in \cite[Section
  8]{ka.hh}. Slightly more generally, we start with an arbitrary
unital monoidal category $\C$, with the product $- \otimes -$ and
the unit object $1 \in \C$.

\begin{defn}
A {\em $\C$-enrichment} of a $2$-category $\E$ is a co-lax
$2$-functor $A:\E^\tau \to B\C^o$ (that is, a lax $2$-functor from
$\E$ to $B\C$). A $\C$-enrichment $A$ of a small category $I$ is a
$\C$-enrichment of its simplicial replacement $\Delta I$, and for
any small category $I'$ and a functor $\gamma:I' \to I$, the {\em
  induced enrichement} $\gamma^*A$ is the composition $A \circ
\Delta\gamma:\Delta I' \to B\C^o$. A {\em morphism} from a
$\C$-enrichement $A_0$ to a $\C$-enrichment $A_1$ is a map $A_0^o
\to A_1^o$ between the opposite functors.
\end{defn}

\begin{exa}
If $I=\ppt$, then a $\C$-enrichment $A:\Delta \to B\C^o$ is the same
thing as an associative unital algebra object in $\C$ --- the object
itself is $A([1]) \in B\C^o_{[1]} \cong \C^o$, and the morphisms
\eqref{fib.fu} for the transition functors \eqref{m.eq} and
\eqref{e.eq} define the product $A([1]) \otimes A([1]) \to A([1])$
and the unit embedding $1 \to A([1])$. Morphisms between
$\C$-enrichments correspond to algebra maps (this is why they are
defined by passing to the opposite functors).
\end{exa}

\begin{exa}\label{A.S.exa}
More generally, if $I=e(S)$ for a set $S$, then a $\C$-enrichment
for $I$ is a small $\C$-enriched category with the set of objects
$S$ in the usual sense: we have a morphism object $A(s,s') \in \C$
for any $s,s' \in S$, the composition maps $A(s,s') \otimes
A(s',s'') \to A(s,s'')$, and the identity maps $1 \to A(s,s)$.
\end{exa}

\begin{exa}\label{A.1.exa}
If $I=[1]$, then a $\C$-enrichment $A$ for $I$ defines algebra
objects $A_{00}$, $A_{11}$ by restricting to $0,1 \in [1]$, and
evaluating $A$ on the tautological object $\id \in \Delta [1]$ of
Example~\ref{d.n.exa} gives an $A_{11}\text{-}A_{00}$-bimodule
$A_{01}$.
\end{exa}

\begin{exa}\label{A.n.exa}
More generally, if $I = [n]$, $n \geq 0$, then a $\C$-enrichment $A$
for $I$ can be visualized as an upper-triangular matrix $(A_{ij})$,
$0 \leq i \leq j \leq n$ of objects in $\C$. For any $0 \leq i \leq
j \leq l \leq n$, we have the product map $A_{ij} \otimes A_{jl} \to
A_{il}$, and these maps are associative and unital in the
appropriate sense.
\end{exa}

For any small category $I$, $\C$-enrichments of $I$ form a full
subcategory in $\Fun_\Delta(\Delta I,B\C^o)^o$, and since
enrichments are functorial with respect to functors $\gamma:I' \to
I$, these categories form a fibration over $\Cat$ whose objects are
pairs $\langle I,A \rangle$. We can consider its restriction to
$\Delta \subset \Cat$, and denote by $\wCat(\C)$ its reduction
\eqref{red.sq}. It is {\em not} our Morita $2$-category, and in
fact, it is not a $2$-category at all: it does not satisfy the Segal
condition. Indeed, already for $[n] = 2$, we have the product map
$A_{01} \otimes A_{12} \to A_{02}$ of Example~\ref{A.n.exa}, but
nothing insures that it factors through an isomorphism $A_{01}
\otimes_{A_{11}} A_{12} \cong A_{02}$ (nor that the tensor product
is even well-defined).

\medskip

To correct for this, consider the functor $\kappa = \rho \circ
\lambda:\Delta \to \Delta$, $[n] \mapsto [n]^<$, and the fibration
$\kappa^*B\C^o$ with the projection $a^*:\kappa^*B\C^o \to B\C^o$
induced by the adjunction map $\Id \to \kappa$. Define a {\em
  module} over a $\C$-enriched small category $\langle I,A \rangle$
as a functor $M:\Delta I \to \kappa^*B\C^o$ over $\Delta$, cartesian
over all initial embeddings $s:[l] \to [n]$ and equipped with an
isomorphism $a^* \circ M \cong A$. Then $\langle I,A
\rangle$-modules form a category, and we denote the opposite
category by $\langle I,A \rangle\amod$. A functor $\gamma:I' \to I$
defines a functor $\gamma^*:\langle I,A \rangle\amod \to \langle
I',\gamma^*A \rangle\amod$, and an object $i \in I_{[0]} \subset
\Delta I$ defines a functor $\langle I,A \rangle\amod \to \C$, $M
\mapsto M(i) \in (\kappa^*B\C^o)^o_{[0]} \cong \C$. Moreover, it has
been proved in \cite[Lemma 8.7]{ka.hh} that if $i \in I$ is an
initial object, then this functor has a left-adjoint $\C \to \langle
I,A \rangle\amod$ sending an object $V \in \C$ to some $\langle I,A
\rangle$-module $V_i$. Modules of the form $V_i$ are called {\em
  representable}.

\begin{defn}\label{A.cyl.def}
A $\C$-enriched category $\langle [n],A \rangle$, $[n] \in \Delta$
is a {\em cylinder} if for any $0 \leq l \leq l' \leq n$, the
functor $t^*:\langle [l'],t^*A \rangle\amod \to \langle [l],t^*A
\rangle\amod$ sends representable modules to representable modules.
\end{defn}

\begin{exa}\label{A.1.cyl.exa}
In the situation of Example~\ref{A.1.exa}, $\langle [1],A \rangle$
is a cylinder iff $A_{01}$ considered as an $A_{11}$-module is of the
form $A_{01} \cong V \otimes A_{11}$ for some $V \in \C$.
\end{exa}

\begin{prop}\label{mor.prop}
The full subcategory $\Mor(\C) \subset \wCat(\C)$ spanned by
cy\-l\-inders is a $2$-category in the sense of
Definition~\ref{2cat.def}.
\end{prop}

\proof{} This is a particular case of \cite[Proposition
  8.20]{ka.hh}, with the following modifications. Firstly, cylinders
of Definition~\ref{A.cyl.def} correspond to iterated polycylinders
of \cite[Definition 8.16]{ka.hh}. Secondly, \cite{ka.hh} encodes
$2$-categories by transpose cofibrations $\C_\perp \to \Delta^o$
rather than fibrations $\C \to \Delta$, see Remark~\ref{cova.rem},
so the notation might be different in places. Thirdly, \cite{ka.hh}
actually deals with a more general situation that combines
Example~\ref{A.S.exa} and Example~\ref{A.1.exa}, and produces the
Morita $2$-category of enriched small categories. To obtain our
$\Mor(\C)$, one has to take the full $2$-subcategory spanned by
algebras (that is, categories with a single object).
\endproof

To see how cylinders of Definition~\ref{A.cyl.def} encode tensor
products, consider the functor $\kappa_\flat:\Delta^< \to \Delta$,
$[n] \mapsto [n]^{<>}$. It sends the initial object $o \in \Delta^<$
to $[1] \in \Delta$, we have $(\kappa_\flat^*B\C)_o \cong
(B\C)_{[1]} \cong \C$, and by Remark~\ref{comma.rem}, the embedding
$\C \cong (\kappa_\flat^*B\C)_o \subset \kappa_\flat^*B\C$ admits a
right-adjoint $\zeta:\kappa_\flat^*B\C \to \C$. Note that for any
$[n] \in \Delta^<$, the unique map $[n] \to [0]$ gives rise to a map
$\kappa_\flat([n]) \to \kappa_\flat([0]) \cong [2]$, and this
promotes the functor $\kappa_\flat$ to a functor
$\kappa_\dg:\Delta^< \to \Delta /^{\Id} [2] \cong \Delta [2]$. Then
for any $\C$-enrichment $A$ of the small category $[2]$, we obtain a
functor $A^\dg = \zeta^o \circ A \circ \kappa_\dg^o:\Delta^{<o} \to
\C$. Its value at the empty ordinal $o \in \Delta^<$ is $A_{02} \in
\C$, and for any $[n] \in \Delta \subset \Delta^<$, we have
\begin{equation}\label{A.dg.eq}
A^\dg([n]) \cong A_{01} \otimes A_{11}^{\otimes n} \otimes A_{12}.
\end{equation}
If $A$ is a cylinder, then $A^\dg$ extends to the category
$\Delta_+^o \supset \Delta^{<o}$, so that $A^\dg \cong
\lambda^{o*}A^\dg_+ \cong \rho^o_!A^\dg_+$ for some
$A^\dg_+:\Delta_+^o \to \C$, where the functor $\rho^o_! \cong
\lambda^{o*}$ is left-adjoint to $\rho^{o*}$. Then since $o \in
\Delta^o_+$ is terminal, we have
$$
A_{02} \cong A^\dg(o) \cong \colim_{\Delta^o_+}A^\dg_+ \cong
\colim_{\Delta^o}\rho^o_!A^\dg_+ \cong \colim_{\Delta^o}A^\dg,
$$
and in particular, the colimit exists. This colimit reduces to the
coequalizer of the two maps $A^\dg(s),A^\dg(t):A_{01} \otimes A_{11}
\otimes A_{12} \to A_{01} \otimes A_{12}$ that encode the two
multiplications, and this is the usual definition of $A_{01}
\otimes_{A_{11}} A_{12}$.

\begin{exa}\label{corr.exa}
For a somewhat tautological but very useful example of a Morita
$2$-category in the sense of Proposition~\ref{mor.prop}, let $\C$ be
a category with finite products, considered as a unital monoidal
category with respect to these products. Then by definition, objects
in $\Mor(\C^o)$ are objects $c \in \C$ equipped with a counital
coalgebra structure, and just as in the symmetric monoidal case,
there is exactly one such for any $c$, with the coproduct $c \to c
\times c$ given by the diagonal map. Thus $\Mor(\C^o)$ has the same
objects as $\C$, and moreover, the morphism category
$\Mor(\C^o)(c,c')$, $c,c' \in \C$ is opposite to the category of
diagrams
\begin{equation}\label{mor.dom}
\begin{CD}
c @<<< \wt{c} @>>> c'
\end{CD}
\end{equation}
in $\C$ such that $\wt{c} \cong x \times c'$ for some $x \in \C$
(this is the cylinder condition). Composition of the diagrams
\eqref{mor.dom} is given by fibered products that exists by virtue
of the cylinder condition. Informally, $\Mor(\C^o)$ is $2$-opposite
to the $2$-category of correspondences in $\C$. Moreover, if we let $\bC
= \C_{\Id}$ be the discrete category underlying $\C$, and let
$e(\C)=e(\bC)$ be the category with the same objects as $\C$ and
exactly one morphism between any two objects, then $\Delta e(\C)
\cong \eps_*\bC$, where $\eps:\ppt \to \Delta$ is the embedding onto
$[0] \in \Delta$, and the identification $\Mor(\C^o)_{[0]} \cong
\bC$ induces a $2$-functor $\Mor(\C^o) \to \Delta e(\C)$. A useful
observation is that it admits a left-adjoint co-lax $2$-functor
\begin{equation}\label{i.C}
i:\Delta e(\C) \to \Mor(\C^o).
\end{equation}
It is identical over $[0]$, and for any $c,c' \in \C$, the unique map
from $c$ to $c'$ goes to the initial object of the category
$\Mor(\C^o)(c,c')$ --- namely, the diagram \eqref{mor.dom} with
$\wt{c} = c \times c'$.
\end{exa}

For a useful generalization of Proposition~\ref{mor.prop}, assume
that the monoidal category $\C$ has all filtered colimits preserved
by the tensor product. Then for any $\C$-enriched category $\langle
I,A \rangle$, $\langle I,A \rangle\amod$ also has filtered colimits,
and one can define ind-representable modules as filtered colimits of
representable ones. Define ind-cylinders by replacing
``representable'' with ``ind-representable'' in
Definition~\ref{A.cyl.def}, and let $\iMor(\C) \subset \wCat(\C)$
be the full subcategory spanned by ind-cylinders. Then
\cite[Proposition 8.20]{ka.hh} also applies to $\iMor(\C)$, so that
it is also a $2$-category. It has the same objects as $\Mor(\C)$ but
more morphisms between them.

\subsection{Derived Morita $2$-category.}\label{dmor.subs}

For a real-life application of Proposition~\ref{mor.prop}, take a
commutative ring $k$, and let $\C = k\amod^{fl}$, the category of
flat $k$-modules. Then objects in $\iMor(k) = \iMor(k\amod^{fl})$
are flat $k$-algebras, and by Example~\ref{A.1.cyl.exa}, morphisms
from $A_0$ to $A_1$ in $\Mor(\C)$ are given by $A_0\text{-}A_1$-bimodules of the
form $V \otimes A_1$, $V \in k\amod^{fl}$. By the standard criterion
of flatness, morphisms in $\iMor(\C)$ are bimodules that are flat
over $A_1$ (``right-flat''). Composition is given by the usual
tensor product over $A_1$.

To generalize this to complexes, we need to generalize slightly
Proposition~\ref{mor.prop}. Assume given unital monoidal categories
$\C$, $\C'$, and a unital monoidal functor $\phi:\C \to \C'$. Say
that $\langle [n],A \rangle \in \wCat(\C)$ is a {\em
  $\phi$-cylinder} resp.\ {\em $\phi$-ind-cylinder} if $B\phi \circ
A$ is a cylinder resp.\ ind-cylinder in the sense of
Definition~\ref{A.cyl.def}, and let $\Mor(\C,\phi),\iMor(\C,\phi)
\subset \wCat(\C)$ be the full subcategories spanned by
$\phi$-cylinders resp.\ $\phi$-ind-cylinders.

\begin{lemma}
Assume that the functor $\phi:\C \to \C'$ is conservative and
creates colimits --- that is, for any small $I$ and functor $c:I \to
\C$, the colimit $\colim_Ic$ exists if and only if so does
$\colim_I\phi(c)$, and in this case, $\colim_I\phi(c) \cong
\phi(\colim_Ic)$. Then $\Mor(\C,\phi)$ is a $2$-category, and
$\wCat(\phi)$ restricts to a $2$-functor $\Mor(\C,\phi) \to
\Mor(\C')$. Moreover, if $\C'$ has filtered colimits, then the same
holds for $\iMor(\C,\phi)$ and $\iMor(\C')$.
\end{lemma}

\proof{} Proposition~\ref{mor.prop} is a part of \cite[Proposition
  8.20]{ka.hh} whose proof relies on \cite[Lemma 8.21]{ka.hh} and
\cite[Lemma 8.22]{ka.hh}. Out of these, Lemma 8.21 holds for
$\phi$-cylinders and $\phi$-ind-cylinders with the same proof, and
so does Lemma 8.22~\thetag{i}. For \thetag{ii} and \thetag{iii} of
Lemma 8.22, notice that since $\phi$ creates colimits, it suffices
to check the existence of $\beta^{\Delta^o}_!$ after applying
$\phi$, and since $\phi$ is also conservative, checking that
$\beta^{\Delta^o}_!\beta^*A \to A$ is an isomorphism can also be
done after applying $\phi$. Finally, the proof of Proposition 8.20
itself for $\phi$-cylinders and $\phi$-ind-cylinders is again
identically the same.
\endproof

For applications, we take $\C = \C_\idot^{fl}(k)$, the category of
termwise-flat chain complexes of $k$-modules, and let $\C' \subset
\C$ be the full subcategory spanned by complexes with zero
differential. Then forgetting the differential gives a conservative
functor $\phi:\C \to \C'$ that creates colimits, so we can form the
$2$-categories $\Mor_\idot(k) = \Mor(\C,\phi)$, $\iMor_\idot(k) =
\iMor(\C,\phi)$. By definition, objects in these $2$-categories are
termwise-flat DG algebras $A_\idot$ over $k$. Morphisms from
$A_\idot$ to $A'_\idot$ in $\Mor_\idot(k)$ resp.\ $\iMor_\idot(k)$
are $A^o_\idot \otimes_k A'_\idot$-modules $M_\idot$ that are
semifree resp.\ semiflat over $A'_\idot$ (that is, become free
resp.\ flat if we forget the differential). We can then consider the
subcategory $\Mor^\flat_\idot(k) \subset \iMor_\idot(k)$ whose
objects are DG algebras that are pointwise-flat and $h$-flat as
complexes of $k$-modules, and whose morphisms are $A^o_\idot
\otimes_k A'_\idot$-modules that are semiflat and $h$-flat over
$A'_\idot$. Since the composition in $\iMor_\idot(k)$ is given by
the usual tensor product, $\Mor^\flat_\idot(k) \subset
\iMor_\idot(k)$ is a $2$-subcategory, and the transition functors of
the fibration $\Mor^\flat_\idot(k) \to \Delta$ send
quasiisomorphisms to quasiisomorphisms. We can then invert those
quasiisomorphisms and obtain the {\em derived Morita $2$-category}
$\DMor(k)$.

\begin{remark}
One usually interprets $\DMor(k)$ as the $2$-category of DG algebras
over $k$ and $k$-linear triangulated functors between their derived
categories but there is a caveat. For any $h$-flat termwise-flat DG
algebras $A_\idot$, $A'_\idot$, $B_\idot$ over $k$, an object $M \in
\DMor(k)(A_\idot,A'_\idot)$ induces a functor
\begin{equation}\label{mor.acts}
- \circ M:\DMor(k)(B_\idot,A_\idot) \to \DMor(k)(B_\idot,A'_\idot),
\end{equation}
and the collection of functors \eqref{mor.acts} determines $M$ (this
is just the Yoneda embedding used in Theorem~\ref{rigi.thm}). If
$B_\idot = k$, then for any $A_\idot \in \DMor(k)_{[0]}$, we have
$\DMor(k)(k,A_\idot) \cong \D(A_\idot)$, so that \eqref{mor.acts} is a
functor $\D(A_\idot) \to \D(A'_\idot)$ between the derived
categories of DG modules. But by itself, this functor is
not enough to recover $M$ --- the correspondence
$\DMor(k)(A_\idot,A'_\idot) \to \Fun(\D(A_\idot),\D(A'_\idot))$ is
neither faithful nor full, and one cannot describe its essential
image in any reasonable way.
\end{remark}

Although we will not need it, it is perhaps useful to describe
$\DMor(k)$ in another way that is somewhat more
homotopy-invariant. Denote by $C_\idot(k)$ the category of all
complexes of $k$-modules, and let $\wCat_\idot(k) \subset
\wCat(C_\idot(k))$ be the full subcategory spanned by objects
$\langle [n],A \rangle$ such that for any map $f:[0] \to [n]$, the
pullback $\langle [0],f^*A \rangle$ lies in $\Mor^\flat_\idot(k)_{[0]}$
(in terms of Example~\ref{A.n.exa}, this means that for any $i \in
[n]$, the DG algebra $A_{ii}$ is termwise-flat and $h$-flat over
$k$).  For any $[n] \in \Delta$, say that a map $\alpha:\langle
[n],A \rangle \to \langle [n],A' \rangle$ in $\wCat_\idot(k)_{[n]}$
is a quasiisomorphism if so is its pullback $f^*\alpha$ with respect
to any map $f:[1] \to [n]$ (in terms of Example~\ref{A.n.exa}, this
means that $\alpha:A_{ij} \to A'_{ij}$ is a quasiisomorphism for any
$0 \leq i \leq j \leq n$). Then note that for any $\langle [2],A
\rangle \in \wCat_\idot(k)_{[2]}$, we have the functor
$A^\dg:\Delta^{<o} \to C_\idot(k)$ of \eqref{A.dg.eq}. This is an
augmented simplicial complex of $k$-modules. We can then take its
normalized chain complex in the usual way, and obtain a bicomplex
$A^\dg_{\idot,\idot}$. Say that $\langle [2],A \rangle$ is a {\em
  quasicylinder} if the sum-total complex $A^\dg_\idot$ of the
bicomplex $A^\dg_{\idot,\idot}$ is acyclic, and more generally, say
that any $\langle [n],A \rangle \in \wCat_\idot(k)_{[n]}$ is a {\em
  quasicylinder} if so is its pullback $\langle [2],f^*A\rangle$ for
any map $f:[2] \to [n]$. By definition, the projection
$\wCat_\idot(k) \to \Delta$ is a fibration whose transition functors
send quasiisomorphisms to quasiisomorphisms and quasicylinders to
quasicylinders.

\begin{lemma}\label{quasi.le}
\begin{enumerate}
\item For any $[n] \in \Delta$, the fiber $\wCat_\idot(k)_{[n]}$ of
  the fibration $\wCat_\idot(k) \to \Delta$ admits a localization
  $\DwCat(k)_{[n]}$ with respect to quasiisomorphisms, so that we
  have a fibration $\DwCat(k) \to \Delta$ with fibers
  $\DwCat(k)_{[n]}$ and transition functors induced by those of
  $\wCat_\idot(k)$.
\item The derived Morita $2$-category $\DMor(k)$ is equivalent to
  the full subcategory in $\DwCat(k)$ spanned by quasicylinders.
\end{enumerate}
\end{lemma}

\proof{} For \thetag{i}, for any DG algebra $A$ termwise-flat and
$h$-flat over $k$, let $T(A)$ be the category of DG algebras $A'$
over $k$ equipped with an additional auxiliary non-negative grading
$A'_\idot$ and an isomorphism $A'_0 \cong A$. Recall that the
category of DG algebras over $k$ carries a model structure whose
fibrations are termwise-surjective maps and whose weak equivalences
are quasiisomorphisms (\cite{H}, \cite[Section 4]{kel}), and note
that the same holds for $T(A)$, with the same proof (to construct
cofibrant replacements, note that the free algebra $T(M,n)$ spanned
by a cofibrant $A^o \otimes_k A$-module $M$ placed in auxiliary
degree $n \geq 1$ is cofibrant by adjunction, and construct a
replacement for some $A' \in T(A)$ by starting with $A$, and adding
the necessary new generators in each auxiliary degree $n$ by
induction on $n$). But for any $[n] \in \Delta$, we have a full
embedding
\begin{equation}\label{A.n.eq}
\wCat_\idot(k)_{[n]} \subset \coprod T(A_{00}
\oplus A_{11} \oplus \dots \oplus A_{nn}), \ A_{00},\dots,A_{nn} \in
\wCat_\idot(k)_{[0]}
\end{equation}
that sends $\langle [n],A \rangle$ to the sum of its components
$A_{ij}$ of Example~\ref{A.n.exa}, with $A_{ij}$ placed in auxiliary
degree $j-i$, and the model structure on its target restricts to a
model structure on its source.

For \thetag{ii}, note that the embedding $\Mor_\idot^\flat(k)
\subset \wCat_\idot(k)$ descends to a functor $\DMor(k) \to
\DwCat(k)$ cartesian over $\Delta$. Let $\DMor(k)' \subset
\DwCat(k)$ be the full subcategory spanned by quasicylinders. Since
the quasicylnder condition is preserved by pullbacks, $\DMor(k)' \to
\Delta$ is a fibration, and we have to check that for any $[n]$, the
functor
\begin{equation}\label{gamma.n}
\gamma_n:\DMor(k)_{[n]} \to \DwCat(k)_{[n]}
\end{equation}
factors through an equivalence $\DMor(k)_{[n]} \cong
\DMor(k)'_{[n]}$. For $n=0$, the statement is tautological, and for
an $n=1$, the quasicylinder condition is empty, so it suffices to
notice that for any termwise-flat $h$-flat DG algebras $A_0$, $A_1$
over $k$, any $A_0^o \otimes_k A_1$-module $M$ is quasiisomorphic to
a module that is semiflat and $h$-flat over $A_1$. However, more is
true: any $M$ is quasiisomorphic to a filtered colimit of iterated
extensions of free modules $A_0^o \otimes_k V \otimes_k A_1$, $V \in
k\amod^{fl}$. Therefore for any $[n] \in \Delta$, any object in
$\Mor^\flat_\idot(k)_{[n]}$ is quasiisomorphic to a filtered colimit
of iterated extensions of cylinders. But for any cylinder $\langle
[2],A \rangle$, the functor $A^\dg$ extends to $\Delta_+^o \supset
\Delta^{<o}$, and this means that the normalized bicomplex
$A^\dg_{\idot,\idot}$ is contractible in one of the directions, so
that its totalization is acyclic. Therefore the embedding
\eqref{gamma.n} sends cylinders to quasicylinders, and since the
quasicylinder condition is closed under extensions and filtered
colimits, the whole \eqref{gamma.n} factors through
$\DMor(k)'_{[n]}$.

Now observe that for any termwise-flat $h$-flat DG algebra $A$ over
$k$, we have the functor $T(A) \to (A^o \otimes_k A)\amod$ sending a
DG algebra $A' \in T(A)$ to its component $A'_1$ of auxiliary degree
$1$, and this functor has a fully faithful left-adjoint $M \mapsto
T(M,1)$. This defines a Quillen adjunction, thus an adjunction of
the localizations, and when we restrict it via the embeddings
\eqref{A.n.eq}, we end up with an adjoint pair consisting of the
embedding \eqref{gamma.n} and a functor $\gamma_n^\dg$ right-adjoint
to it. Therefore firstly, \eqref{gamma.n} is fully faithful, and
secondly, its essential image consists of objects $A$ such that the
adjunction map $a:\gamma_n\gamma_n^\dg A \to A$ is an isomorphism in
$\DwCat(k)_{[n]}$.

To finish the proof, it suffices to show that this holds for any
quasicylinder, and since $\gamma^\dg_n(a)$ is an isomorphism by
adjunction, it suffices to show that $\gamma^\dg_n$ is conservative
on $\DMor(k)'_{[n]}$. For $n=2$, this is clear: if a morphism
$\alpha:A \to A'$ induces quasiisomorphisms $A_{01} \to A_{01}'$,
$A_{11} \to A_{11}'$, $A_{12} \to A_{12}'$, then it induces
quasiisomorphisms on all the terms in \eqref{A.dg.eq}, and since the
total complexes are acyclic, the component $A_{02} \to A'_{02}$ must
also be a quasiisomorphism. The general case then follows by
induction on $n$.
\endproof

\begin{remark}\label{quasi.rem}
Note that the quasicylinder condition is manifestly
symmetric. Therefore Lemma~\ref{quasi.le} immediately implies that
we have an equivalence
\begin{equation}\label{mor.io}
\iota:\DMor(k) \cong \DMor(k)^\iota
\end{equation}
sending a DG algebra $A_\idot$ to the opposite DG algebra
$A^o_\idot$.
\end{remark}

\subsection{Fourier-Mukai $2$-category.}

Let us now combine Example~\ref{corr.exa} with the technology of
Subsection~\ref{dmor.subs}, and apply it to the categories of
Example~\ref{sch.exa}, with the monoidal structures of
Subsection~\ref{comm.subs}. Assume that $k$ is Noetherian. The
category $\Sch(k)$ has finite products, so that
Example~\ref{corr.exa} applies directly, and shows that the Morita
$2$-category $\Mor(\Sch(k)^o)$ is $2$-opposite to the $2$-category
of correspondences in $\Sch(k)$. Then as in
Subsection~\ref{dmor.subs}, we let $\C \subset \QCoh_\idot(\Sch(k))$
be the full subcategory spanned by pairs $\langle X,\F_\idot
\rangle$ such that the complex $\F_\idot$ of quasicoherent sheaves
of $\calo_X$-modules is termwise-flat, we consider the forgetfull
functor $\phi:\C \to \C'$, where $\C' \subset \C$ consists of
complexes with trivial differential, and we construct the Morita
$2$-category $\Mor(\C,\phi)$. Moreover, while $\C$ does not
necessarily have all filtered colimits, the fibers of the
cofibration $\C \to \Sch(k)^o$ do; we add those filtered colimits of
$\phi$-cylinders to $\Mor(\C,\phi)$ and by abuse of notation, denote
the resulting $2$-category by $\iMor(\C,\phi)$. Then by definition,
objects in $\iMor(\C,\phi)$ are pairs $\langle X,\A_\idot \rangle$
of a scheme $X \in \Sch(k)$ and a DG algebra $\A_\idot$ of flat
quasicoherent sheaves on $X$, and morphisms from $\langle X,\A_\idot
\rangle$ to $\langle X',\A'_\idot \rangle$ are pairs $\langle
Z,\M_\idot \rangle$ of a correspondence $Z \to X \times X'$ in
$\Sch(k)$ abstractly isomorphic over $X'$ to $Y \times X'$ for some
$Y \in \Sch(k)$, and a sheaf of DG $\pi^*\A_\idot^o
\otimes_{\calo_Z} {\pi'}^*\A'_\idot$-modules on $Z$, semiflat over
${\pi'}^*\A'_\idot$, where $\pi:Z \to X$, $\pi':Z \to Z$ are the
projections. If we have another object $\langle X'',\A''_\idot
\rangle$ and a morphism $\langle Z',\M'_\idot \rangle$ from $\langle
X',\A'_\idot \rangle$ to $\langle X'',\A''_\idot \rangle$, then the
composition is given by
\begin{equation}\label{bmor.compo}
\langle Z,\M_\idot \rangle \circ \langle Z',\M'_\idot \rangle \cong
\langle Z \times_{X'} Z',\nu^*\M_\idot \otimes_{\eta^*\A'_\idot}
{\nu'}^*\M'_\idot \rangle,
\end{equation}
where $\nu:Z \times_{X'} Z' \to Z$, $\nu':Z \times_{X'} Z' \to Z'$,
$\eta = \pi'\circ \nu:Z \times_{X'} Z' \to X'$ are again the
projections.

We then observe that by virtue of \eqref{bmor.compo}, objects
$\langle X,\A_\idot \rangle$ such that $\A_\idot$ is $h$-flat over
$X$, and morphisms $\langle Z,\M_\idot \rangle$ such that $\M_\idot$
is $h$-flat over $Z$ form a $2$-subcategory
$\Mor_\idot^\flat(\Sch(k)) \subset \iMor(\C,\phi)$. We then have a
pair of $2$-opposite $2$-functors
\begin{equation}\label{mor.ff}
\begin{aligned}
\psi:&\Mor_\idot^\flat(\Sch(k)) \to \Mor(\Sch(k)^o),\\
\psi^\tau:&\Mor_\idot^\flat(\Sch(k))^\tau \to \Mor(\Sch(k)^o)^\tau,
\end{aligned}
\end{equation}
where $\psi$ is induced by the forgetful functor
$\QCoh_\idot(\Sch(k)) \to \Sch(k)^o$.

\begin{lemma}\label{psi.t.le}
The functor $\psi^\tau$ of \eqref{mor.ff} is a fibration.
\end{lemma}

\proof{} Use Lemma~\ref{fib.le}. The conditions \thetag{i},
\thetag{ii} hold by definition. For \thetag{iii}, note that
$\psi^\tau$ is an identity over $[0] \in \Delta$, and its fiber over
some $Z \to X \times X'$ in $\Mor(\Sch(k)^o)^\tau_{[1]}$ is given by
$$
\Mor_\idot^\flat(\Sch(k))^\tau_Z \cong \coprod_{\A_\idot,\A'_\idot}
\QCoh^\flat_\idot(Z,\pi^*\A_\idot^o \otimes_{\calo_Z}
{\pi'}^*\A'_\idot)^o,
$$
where the disjoint union is over all pairs of $h$-flat termwise-flat
DG algebras $\A_\idot \in \QCoh_\idot(X)$, $\A'_\idot \in
\QCoh_\idot(X')$, and $\QCoh^\flat_\idot(-,-)$ stands for the
category of DG modules $h$-flat over $\calo_Z$ and semiflat over
${\pi'}^*\A'_\idot$. Then $\psi^\tau([1])$ is obviously a fibration,
with the transition functor
\begin{equation}\label{Z.M}
f^*(\langle Z,\M_\idot \rangle) = \langle W,f^*\M_\idot \rangle
\end{equation}
for any map $f:W \to Z$ of correspondences from $X$ and $X'$, and
by the Segal condition, $\psi^\tau([n])$ is a fibration for any
$[n]$. Then again by the Segal condition, it suffices to check
\thetag{iv} for the map $m:[1] \to [2]$ in $\Delta$ of \eqref{m.eq},
and in this case, it immediately follows from \eqref{Z.M} and
\eqref{bmor.compo}.
\endproof

Now note that the transition functors \eqref{Z.M} of the fibration
$\psi^\tau$ provided by Lemma~\ref{psi.t.le} obviously preserve
quasiisomorphisms. Therefore we can again localize with respect to
quasiisomorphims and obtain a category $\E$ fibered over
$\Mor(\Sch(k)^o)^\tau$, with the fiber over some $Z \to X \times X'$
given by the disjoint union of the derived categories
$\D(Z,\pi^*\A_\idot \otimes_{\calo_X} {\pi'}^*\A'_\idot)$. The
category $\E$ satisfies the Segal condition, and defines a
$2$-category in the sense of Definition~\ref{2cat.def}. We denote
the $2$-opposite $2$-category by $\DMor(\Sch(k))$ and call it the
{\em big derived Morita $2$-category} over $k$. The functor $\psi$
of \eqref{mor.ff} then induces a $2$-functor
\begin{equation}\label{psi.eq}
\psi:\DMor(\Sch(k)) \to \Mor(\Sch(k)^o).
\end{equation}
Objects of $\DMor(\Sch(k))$ are again pairs $\langle X,\A_\idot
\rangle$, and morphisms are pairs $\langle Z,\M \rangle$, $\M \in
\D(Z,\pi^*\A_\idot^o \otimes_{\calo_Z} {\pi'}^*\A'_\idot)$, where $Z
\cong Y \times X'$ for some $Y$. Composition is given by the derived
version of \eqref{bmor.compo}.

\begin{lemma}\label{psi.le}
The $2$-functor \eqref{psi.eq} is a fibration.
\end{lemma}

\proof{} As in Lemma~\ref{psi.t.le}, the conditions \thetag{i},
\thetag{ii} of Lemma~\ref{fib.le} are satisfied by definition, and
it suffices to check \thetag{iii} over $[1] \in \Delta$ where it
amounts to observing that when we pass to the derived categories,
the transition functor \eqref{Z.M} admits a right-adjoint
$f_*$. Moreover, it suffices to check \thetag{iv} for the map $m$ of
\eqref{m.eq}. To do this, assume given objects $\langle
X,\A_\idot\rangle$, $\langle X',\A'_\idot\rangle$, $\langle
X'',\A''_\idot\rangle$ in $\DMor(\Sch(k))_{[0]}$ and morphisms
$$
\begin{aligned}
f:\langle W,\M_\idot \rangle \to \langle Z,f_*\M_\idot \rangle
&\in \DMor(\Sch(k))(\langle X,\A_\idot \rangle,\langle X',\A'_\idot
\rangle),\\
f':\langle W',\M'_\idot \rangle \to \langle Z',f'_*\M'_\idot \rangle 
&\in \DMor(\Sch(k))(\langle X',\A'_\idot \rangle,\langle X'',\A''_\idot
\rangle),
\end{aligned}
$$
both cartesian over $\Mor(\Sch(k)^o)_{[1]}$. Let $\nu$, $\nu'$ and
$\eta$ have the same meaning as in \eqref{bmor.compo}, and denote by
$\eps:W \times_{X'} W' \to W$, $\eps':W \times_{X'} W' \to W'$ and
$\chi = \nu \circ (f \times f'):W \times_{X'} W \to X'$ the
corresponding projections for $W$ and for $W'$. By
\eqref{bmor.compo}, we then have to check that the base change map
\begin{equation}\label{bc.1}
\nu^*f_*\M_\idot \lotimes_{\eta^*\A'_\idot} {\nu'}^*f'_*\M'_\idot \to
(f \times f')_*\eps^*\M_\idot \lotimes_{\chi^*\A'_\idot}
   {\eps'}^*\M'_\idot
\end{equation}
is an isomorphism. The pushforward functors are triangulated, and
since we are working with schemes of finite type, they also commute
with arbitrary sums. Therefore we may assume that $\M_\idot$ is free
over $\chi^*\A'_\idot$ --- that is, we have $\M_\idot \cong \F_\idot
\lotimes_{\calo_W} {\chi'}^*\A'_\idot$ for some $\F_\idot \in
\D(W)$. Then by the projection formula, we have $f_*\M_\idot \cong
f_*\F_\idot \lotimes_{\calo_Z} {\nu'}^*\A'_\idot$, and \eqref{bc.1}
can be rewritten as
\begin{equation}\label{bc.2}
\nu^*f_*\F_\idot \lotimes_{\calo_{Z \times_{X'} Z}}
   {\nu'}^*f'_*\M'_\idot \to (f \times f')_*\eps^*\M_\idot
   \lotimes_{\calo_{W \times_{X'} W'}} {\eps'}^*\M'_\idot.
\end{equation}
Moreover, we have $(f \times f') = (f \times \id) \circ (\id \times
f')$, so it suffices to consider separately the cases $f=\id$ and
$f'=\id$. In the first case, we may again apply the projection
formula to reduce to the case $\F_\idot \cong \calo_W$, and then
\eqref{bc.2} becomes the base change map $\nu^*f_*\M'_\idot \to (f
\times \id)_*\eps^*\M'_\idot$ associated to the cartesian square
\begin{equation}\label{bc.sq.1}
\begin{CD}
Z \times_{X'} W' @>{f \times \id}>> Z \times_{X'} Z'\\
@V{\eps}VV @VV{\nu}V\\
W' @>{f}>> Z',
\end{CD}
\end{equation}
and since both $W$ and $Z$ are cylinders, the vertical maps are
flat, so we are done by flat base change. In the second case, we
reduce to the case $\M'_\idot \cong \calo_{W'}$, and we note that
while the vertical maps in the corresponding version of the
cartesian square \eqref{bc.sq.1} are not flat, they become flat if
we replace the products over $X'$ with the absolute products. Then
it suffices to check that the leftmost square in the diagram
$$
\begin{CD}
W \times_{X'} Z' @>{f \times \id}>> Z \times_{X'} Z' @>{\nu'}>> Z'\\
@VVV @VVV @VVV\\
W \times Z' @>{f \times \id}>> Z \times Z' @>>> X' \times Z'
\end{CD}
$$
satisfies base change. But both squares are cartesian, and both maps
$W \times Z',Z \times Z' \to X' \times Z'$ are flat, so we are again
done by flat base change.
\endproof

Now consider the co-lax $2$-functor \eqref{i.C} for $2$-the category
$\Mor(\Sch(k)^o)$, and note that by virtue of Lemma~\ref{psi.le},
$\FMbig(k) = i^*\DMor(\Sch(k))$ comes equipped with a fibration
$\FMbig(k) \to \Delta e(\Sch(k)) \to \Delta$ that turns it into a
$2$-category in the sense of Definition~\ref{2cat.def}. We call it
the {\em big Fourier-Mukai $2$-category} over $k$. Objects are again
pairs $\langle X,\A_\idot \rangle$ of a scheme $X$ and an $h$-flat
termwise-flat DG algebra $\A_\idot$, and the categories of morphisms
are
\begin{equation}\label{FM.eq}
\FMbig(k)(\langle X,\A_\idot \rangle,\langle X',\A'_\idot \rangle)
\cong \D(X \times
X',\pi^*\A^o_\idot \otimes_{\calo_{X \times X'}}
{\pi'}^*\A'_\idot),
\end{equation}
where $\pi:X \times X' \to X$, $\pi':X \times X' \to X'$ are the two
projections. The composition $\M_\idot \circ \M'_\idot$ of morphisms
$\M_\idot \in \FMbig(\langle X,\A_\idot\rangle,\langle
X',\A'_\idot\rangle)$ and $\M'_\idot \in \FMbig(\langle
X',\A'_\idot\rangle,\langle X'',\A''_\idot\rangle)$ is given by
\begin{equation}\label{FM.compo}
\M_\idot \circ \M'_\idot \cong \rho_*(\nu^*\M_\idot
\lotimes_{\eta^*\A'_\idot} {\nu'}^*\M'_\idot),
\end{equation}
where $\nu$, $\nu'$ and $\eta$ are as in \eqref{bmor.compo}, and
$\rho:X \times X' \times X'' \to X \times X''$ is the projection. We
note that we have an obvious full embedding $\DMor(k) \subset
\FMbig(k)$ that identifies $\DMor(k)$ with the $2$-subcategory in
$\FMbig(k)$ spanned by pairs $\langle X,\A_\idot \rangle$ with $X =
\ppt = \Spec k$. We also have the {\em Fourier-Mukai $2$-category}
defined as the full $2$-subcategory $\FM(k) \subset \FMbig(k)$
spanned by pairs $\langle X,\calo_X \rangle$. This is the same
Fourier-Mukai $2$-category as the one discussed in \cite[Subsection
  2.6]{AL2}. We note that as in Remark~\ref{quasi.rem},
\eqref{FM.eq} is manifestly symmetric, so that we have an
equivalence
\begin{equation}\label{fm.io}
\iota:\FMbig(k) \cong \FMbig(k)^\iota
\end{equation}
sending $\langle X,\A_\idot \rangle$ to $X,\A_\idot^o \rangle$.

\subsection{Equivalences.}

According to our definition, two quasiisomorphic DG algebras over
$k$ are different as objects in the derived Morita $2$-category
$\DMor(k)$. However, they are equivalent in the general
$2$-categorical sense. More precisely, assume given a map $f:A_\idot
\to B_\idot$ between $h$-flat termwise-flat DG algebras $A_\idot$,
$B_\idot$. Define the {\em graph} $\gr(f)$ of the map $f$ as the
object
\begin{equation}\label{gr.f}
\gr(f) = A_\idot \lotimes_{A_\idot^o \otimes_k A} A^o_\idot
\otimes_k B_\idot \in \D(A_\idot^o \otimes_k B_\idot) \cong
\DMor(k)(A_\idot,B_\idot),
\end{equation}
where $A_\idot$ on the left is the diagonal bimodule, and
$A^o_\idot \otimes_k A_\idot$ acts on $A^o_\idot \otimes_k B_\idot$
via the map $\id \otimes f$. Then we have $\gr(f \circ g) \cong
\gr(f) \circ \gr(g)$ for any composable pair of maps $f$, $g$, and
if $f$ is a quasiisomorphism, we have
\begin{equation}\label{gr.f.f}
\begin{aligned}
\gr(f) \circ \iota(\gr(f^o)) &\cong B_\idot \in \D(B^o_\idot
\otimes_k B_\idot),\\
\iota(\gr(f^o)) \circ \gr(f) &\cong A_\idot \in \D(A^o_\idot \otimes_k
A_\idot),
\end{aligned}
\end{equation}
where $\iota$ is the equivalence \eqref{mor.io} of
Remark~\ref{quasi.rem}. Therefore $\gr(f)$ and $\iota(\gr(f^o))$ are
a pair of mutually inverse equivalences. Not all the equivalences in
$\DMor(k)$ are of this form, but those that are can be easily
characterized. Namely, both $A_\idot$ and $B_\idot$ define algebra
objects $A$, $B$ in the monoidal category $\D(k)$, and any object
$M_\idot \in \D(A^o_\idot \otimes_k B_\idot)$ defines an object $M
\in \D(k)$ that is a left module over $B$ and a right module over
$A$. In particular, any map $e:k \to M$ defines action maps
\begin{equation}\label{LR}
L(k):A = A \otimes k \to M, \qquad R(k):B = B \otimes k \to M
\end{equation}
in the derived category $\D(k)$.

\begin{lemma}\label{mor.eq.le}
An object $M_\idot \in \D(A^o_\idot \otimes_k B_\idot)$ is of the
form $\gr(f)$ for a quasiisomorphism $f:A_\idot \to B_\idot$ if and
only if for some map $e:k \to M$ in $\D(k)$, both maps \eqref{LR}
are isomorphisms.
\end{lemma}

\proof{} The ``only if'' part is immediate from \eqref{gr.f} (since
it is obvious for the diagonal bimodule). For the converse,
represent $M$ by an $h$-projective $A_\idot^o \otimes_k
B_\idot$-module $M_\idot$, lift $e$ to a map of complexes $e:k \to
M_\idot$, and note that it induces a quasiisomorphism $B_\idot =
B_\idot \otimes_k k \to M_\idot$ of $B_\idot$-modules. Then the
action of $A_\idot$ and evaluation at $e(1) \in M_0$ give maps of
complexes
\begin{equation}\label{act.dia}
\begin{CD}
A^o_\idot @>{f}>> \End^\hdot_{B_\idot^\hdot}(M_\idot)^o
@>{\ev}>> M_\idot,
\end{CD}
\end{equation}
the map $f$ is a DG algebra map, and since $M_\idot$ is
$h$-projective, the target of the map $f$ is quasiisomorphic to
$B_\idot$. Then $M \cong \gr(f)$, and since $\ev$ and $\ev
\circ f$ in \eqref{act.dia} are quasiisomorphisms by assumption, so
is the map $f$.
\endproof

Now let us turn to the Fourier-Mukai $2$-category
$\FMbig(k)$. Recall that an object $\E \in \C$ in a unital monoidal
category $\C$ is {\em (left-)dualizable} if it is reflexive as a
morphism in the $2$-category $B\C$; in this case, we have the
adjoint morphism $\E^\vee \in (B\C)_{[1]} \cong \C$ and the
adjunction maps $\E \circ \E^\vee \to 1$, $1 \to \E^\vee \circ
\E$. Recall also that for any $X \in \Sch(k)$, an object $\E \in
\D(X)$ is dualizable if and only if it is perfect (that is,
$\Hom(\E,-)$ commutes with arbitrary sums). Finally, recall that $\E
\in \D(X)$ is a {\em generator} if $\Hom(\E,\F) = 0$ implies $\F=0$,
and recall that by \cite{BvB}, any $X \in \Sch(X)$ admits a perfect
generator $\E \in \D(X)$.

For any $\E \in \D(X)$, we can represent $\E$ by an $h$-injective
complex $\E_\idot$ of quasicoherent sheaves on $X$, and then choose
a cofibrant $h$-flat termwise-flat replacement $A(\E)_\idot \to
\End^\hdot(\E_\idot)$ of the DG algebra of endomorphisms of the
complex $\E_\idot$. Then $A(E)_\idot$ acts on $\E_\idot$, so we can
promote $E$ to an object
\begin{equation}\label{K.E}
\K(E) \in \D(X,A(\E)_\idot \otimes_k \calo_X) = \FMbig(k)(\langle
X,\calo_X \rangle, \langle \ppt,A(\E)_\idot\rangle).
\end{equation}
If $\E$ is perfect, thus dualizable, we have the dual object
$\E^\vee$, and we can do the same procedure to obtain an $h$-flat
termwise-flat DG algebra $A(E^\vee)_\idot$ and a object
\begin{equation}\label{K.E.vee}
\K(\E^\vee) \in \D(X,A(\E^\vee)_\idot \otimes_k \calo_X) =
\DMor(k)(\langle \ppt,A(\E^\vee)^o_\idot \rangle, \langle X,\calo_X
\rangle)
\end{equation}
where we use the equivalence \eqref{fm.io} to interpret
$\K(\E^\vee)$ as a morphism to and not from $\langle X,\calo_X
\rangle \cong \langle X,\calo_X^o \rangle$.

\begin{lemma}\label{fm.eq.le}
For any $X \in \Sch(h)$ and perfect $E \in \D(X)$, the objects
\eqref{K.E} and \eqref{K.E.vee} satisfy $\K(\E) \circ \K(\E^\vee)
\cong \id$. If $\E \in \D(X)$ is a generator, then also $\K(\E^\vee)
\circ \K(\E) \cong \id$.
\end{lemma}

\proof{} For the first claim, note by \eqref{FM.compo}, the object
$M \in \D(k)$ underlying $\K(\E) \circ \K(\E^\vee) \in
\D(A(\E^\vee)_\idot \otimes_k A(\E)_\idot)$ is given by
$$
M \cong H^\hdot(X,\E \lotimes_{\calo_X} \E^\vee),
$$
the algebras $\End^\hdot(\E)$ resp.\ $\End^\hdot(\E^\vee)$ act on
$M$ via $\E$ resp.\ $\E^\vee$, and then the adjunction map $\calo_X
\to \E \lotimes_{\calo_X} \E^\vee$ induces a map $e:k \to M$ such
that both $L(k):\End^\hdot(\E) \to M$ and $R(k):\End^\hdot(\E^\vee)
\to M$ are isomorphisms. We are then done by Lemma~\ref{mor.eq.le},
and as a bonus, we obtain a quasiisomorphism between
$A(\E^\vee)_\idot^o$ and $A(\E)_\idot$.

For the second claim, note that composition with $\K(E)$ gives a
functor
$$
\begin{aligned}
L(\E):\D(A(\E)^o_\idot) \cong &\FMbig(\langle \ppt,A(\E)_\idot
\rangle,\langle \ppt,k \rangle) \longrightarrow\\ 
&\longrightarrow \FMbig(\langle X,\calo_X \rangle, \langle \ppt,k
\rangle) \cong \D(X),
\end{aligned}
$$
and this functor has a right-adjoint $R(\E):\D(X) \to
\D(A(\E)^o_\idot)$ sending some $\F_\idot$ to
$\Hom^\hdot(\E_\idot,\F_\idot)$ with the natural $A(\E)_\idot$-module
structure. If $\E$ is a generator, then $L(\E)$ and $R(\E)$ are an
inverse pair of equivalences. But in this case, it is well-known
that $\E^\vee$ is also a generator, and $E \boxtimes_k \E^\vee$ is a
generator for $\D(X \times X)$. We then have a adjoint pair of
equivalences
$$
\begin{aligned}
L(\E \boxtimes \E^\vee):&\D(A(\E)^o_\idot \otimes_k
A(\E^\vee)^o_\idot) \to \D(X \times X),\\
R(\E \boxtimes \E^\vee):&\D(X \times X) \to \D(A(\E)^o_\idot \otimes
A(\E^\vee)^o_\idot).
\end{aligned}
$$
If we identify $A(\E^\vee)^o_\idot \cong A(\E)_\idot$, then $L(\E
\boxtimes \E^\vee)$ sends the diagonal bimodule $A(\E)_\idot$ to the
composition $\K(\E^\vee) \circ \K(\E) \in \D(X \times X)$, while the
adjoint and inverse equivalence $R(\E \boxtimes \E^\vee)$ sends the
structure sheaf $\calo_X$ of the diagonal $X \subset X \times X$ to
the diagonal bimodule $A(\E)_\idot$. Therefore $\K(\E) \circ \K(\E^\vee)
\cong \calo_X$, and this is exactly the identity morphism in
$\FMbig(k)$.
\endproof

Lemma~\ref{fm.eq.le} has a useful corollary. Say that a {\em
  framing} for the $2$-category $\FM(k)$ is a choice of a generator
$E \in \D(X)$, its $h$-injective representative $\E_\idot$, and a
cofibrant $h$-flat termwise-flat replacement $A(\E)_\idot \cong
\End^\hdot(\E_\idot)$ for any $X \in \Sch(k)$.

\begin{corr}\label{fm.corr}
Any framing for $\FM(k)$ defines a fully faithful $2$-functor
$\FM(k) \to \DMor(k)$ sending $\langle X,\calo_X \rangle$ to
$A(\E)_\idot$, and two different framings define equivalent
$2$-functors.
\end{corr}

\proof{} For any framing, we can take the equivalences of
Lemma~\ref{fm.eq.le} and plug them into Proposition~\ref{tw.prop} as
in Example~\ref{equi.exa}; the resulting twisting functor
$\Theta:\FM(k) \to \FMbig(k)$ factors through $\DMor(k)$ and gives
the required fully faithful $2$-functor. If we have two framings, we
can combine them to obtain a collection of equivalences for the
$2$-category $\FM(k) \times^2 \eq$, and the twisting functor then
gives an equivalence in the sense of \eqref{ga.adj.eq}.
\endproof

{\small\noindent
Affiliations:
\begin{enumerate}
\renewcommand{\labelenumi}{\arabic{enumi}.}
\item Steklov Mathematics Institute, Algebraic Geometry Section
  (main affiliation).
\item Laboratory of Algebraic Geometry, National Research University
Higher\\ School of Economics.
\end{enumerate}}

{\small\noindent
{\em E-mail address\/}: {\tt kaledin@mi-ras.ru}
}
\end{document}